\documentclass[12pt]{article}
\usepackage{geometry}
\usepackage[T1]{fontenc}     
\usepackage{lmodern}         

\usepackage{amssymb}
\newcounter{bla}

\newtheorem{proposition}{Proposition}[section]  

\usepackage[font=normalsize]{caption}
\usepackage{subcaption} 
\usepackage{amsmath}
\usepackage[numbers]{natbib}  
\usepackage[pdftex]{graphicx}
\usepackage{calc}
\usepackage{multicol}
\usepackage{multirow}
\usepackage{mathtools}
\usepackage{siunitx}
\usepackage{booktabs}
\usepackage{rotating}
\usepackage[colorlinks]{hyperref}
\hypersetup{
	colorlinks,
	citecolor=magenta,
	linkcolor=cyan,
	filecolor=blue,      
	urlcolor=red,
	pdftitle={MLMC covariance},
	pdfpagemode=FullScreen,
}
\usepackage{etoolbox}
\usepackage{tikz}
\usepackage{url}
\usepackage{caption}
\usepackage{cleveref}
\usepackage[export]{adjustbox}
\allowdisplaybreaks

\usepackage{color}
\definecolor{mygreen}{rgb}{0.45, 0.65, 0.0}
\newcommand{\rvsi}[1]{\textcolor{black}{#1}}
\definecolor{myblue}{rgb}{0.1, 0.2, 0.9}
\newcommand{\rvsn}[1]{\textcolor{black}{#1}}
\definecolor{myred}{rgb}{0.8, 0.1, 0.0}

\newcommand{\rvso}[1]{\textcolor{black}{#1}}
\usepackage{pgfplots}
\usepgfplotslibrary{groupplots,dateplot}
\usetikzlibrary{patterns,shapes.arrows}
\pgfplotsset{compat=newest}

\usepackage[sfdefault=cmbr,OMLmathsans]{isomath}  

\newcommand{\vek}[1]{\mathchoice{\displaystyle\boldsymbol#1}
	{\textstyle\boldsymbol#1}{\scriptstyle\boldsymbol#1}
	{\scriptscriptstyle\boldsymbol#1}}

\newcommand{\mrm}{\mathrm}     

\usepackage[nopostdot,nonumberlist,acronym,toc,nogroupskip]{glossaries}
\setglossarystyle{long}

\glsdisablehyper 
\newacronym{bmd}{BMD}{bone mineral density}
\newacronym{qct}{QCT}{quantitative computed tomography}
\newacronym{ct}{CT}{computed tomography}
\newacronym{hu}{HU}{Hounsfield units}
\newacronym{dxa}{DXA}{dual energy X-ray absorptiometry}
\newacronym{fea}{FEA}{finite element analysis}
\newacronym{fem}{FEM}{finite element method}
\newacronym{mc}{MC}{Monte Carlo method}
\newacronym{mpp}{MPP}{most probable point}
\newacronym{amv}{AMV}{advanced mean-value method}
\newacronym{form}{FORM}{first order reliability method}
\newacronym{sorm}{SORM}{second order reliability method}
\newacronym{rsms}{RSMs}{response surface methods}
\newacronym{pce}{PCE}{polynomial chaos expansion}
\newacronym{bc}{BC}{boundary condition}
\newacronym{pdf}{PDF}{probability density function}
\newacronym{mse}{MSE}{mean square error}
\newacronym{nmse}{NMSE}{normalized mean square error}
\newacronym{rmse}{RMSE}{root mean square error}
\newacronym{lhs}{LHS}{Latin hypercube sampling}
\newacronym{ism}{ISM}{importance sampling method}
\newacronym{kle}{KLE}{Karhunen-Lo\`{e}ve expansion}
\newacronym{spde}{SPDE}{stochastic partial differential equation}
\newacronym{pde}{PDE}{partial differential equation}
\newacronym{mlmc}{MLMC}{multilevel Monte Carlo method}
\newacronym{smlmc}{SMLMC}{scale invariant multilevel Monte Carlo method}
\newacronym{dof}{DOF}{degrees of freedom}
\newacronym{iso}{iso}{isotropy}
\newacronym{ortho}{ortho}{orthotropy}
\newacronym{scl}{scl}{scaling}
\newacronym{dir}{dir}{direction}
\newacronym{uq}{UQ}{uncertanity quantification}
\newacronym{pl}{PL}{Plaston}
\newacronym{nt}{NT}{nodal temperature}
\newacronym{thfl}{THFL}{total heat flux}
\newacronym{nhfl}{NHFL}{normalized heat flux}
\newacronym{td}{TD}{total displacement}
\newacronym{rmt}{RMT}{random matrix theory}
\newacronym{trf}{TRF}{tensor-valued random field}
\newacronym{spd}{SPD}{symmetric and positive definite}
\usepackage[intoc, english]{nomencl}
\renewcommand{\nompreamble}{The important notations used in this thesis are listed here, in the sequence in which they were introduced in the chapters:}
\usepackage{etoolbox}
\renewcommand\nomgroup[1]{%
  \item[\bfseries
  \ifstrequal{#1}{A}{Chapter 2}{%
  \ifstrequal{#1}{B}{Chapter 3}{%
  \ifstrequal{#1}{C}{Chapter 4}{%
  \ifstrequal{#1}{D}{Chapter 5}{}}}}%
]}
\makenomenclature

\newcommand{\spatialdomain}{\mathcal{G}}

\newcommand{\Euclideanspace}{\mathbb{R}^d}

\newcommand{\solspace}{\mathcal{U}}

\newcommand{\elesize}{h}

\newcommand{\convorder}{\alpha}
\newcommand{\convbeta}{\beta}
\newcommand{\convgamma}{\gamma}

\newcommand{\samplespace}{\varOmega}
\newcommand{\event}{\omega}

\newcommand{\sigmaalgebra}{\mathfrak{F}}
\newcommand{\probability}{\mathbb{P}}

\newcommand{\samplesize}{N}

\newcommand{\normaldist}{\mathcal{N}}

\newcommand{\meshlevel}{l}
\newcommand{\ml}{\meshlevel}
\newcommand{\numofmesh}{\LL}

\newcommand{\optSam}{N_l}

\newcommand{\Lagrange}{\tau}

\newcommand{\SymMatSpace}{\mrm{Sym}}
\newcommand{\PosSymMatSpace}{\SymMatSpace^+}

\newcommand{\nodes}{M}

\newcommand{\uu}{u} 
\newcommand{\h}{\mathrm{h}}
\newcommand{\N}{\samplesize}
\newcommand{\uh}{\uu^h}

\newcommand{\U}{\mathcal{U}}
\newcommand{\mumu}{\mu} 
\newcommand{\muMC}{\widehat{\mumu}^{\mathrm{MC}}}
\newcommand{\htwoMC}{\widehat{\h}_2^{\mathrm{MC}}}

\newcommand{\Nl}{\N_\meshlevel}
\newcommand{\LL}{L} 

\newcommand{\uhL}{\uu^{h_\LL}}

\newcommand{\Zl}{{Z}_{\meshlevel}}
\newcommand{\Yl}{{Y}_{\meshlevel}}

\newcommand{\epssq}{\epsilon^2}
\newcommand{\epssqtwo}{\epssq/2}

\newcommand{\PosRealStrict}{\mathbb{R}^{+}}

\newcommand{\heatsource}{f}
\newcommand{\temperature}{\uu}

\newcommand{\var}{\mathbb{V}\mrm{ar}}
\newcommand{\varMC}{\widehat{\mathbb{V}}\mrm{ar}^{\mrm{MC}}}

\newcommand{\cov}{\mrm{cov}}
\newcommand{\covMC}{\mrm{c}\widehat{\mrm{o}}\mrm{v}^{\mrm{MC}}}
\newcommand{\covML}{\mrm{c}\widehat{\mrm{o}}\mrm{v}^{\mrm{ML}}}
\newcommand{\aone}{{a}_1}
\newcommand{\atwo}{{a}_2}
\newcommand{\athree}{{a}_3}
\newcommand{\uhij}{(u^h_i, u^h_j)}
\newcommand{\uij}{(u_i, u_j)}
\newcommand{\s}{s}
\newcommand{\honeoneMC}{\widehat{\h}_{1,1}^{\mathrm{MC}}}
\newcommand{\honeoneML}{\widehat{\h}_{1,1}^{\mathrm{ML}}}
\newcommand{\hrMC}{\widehat{\h}_\p^{\mathrm{MC}}}
\newcommand{\htwotwoMC}{\widehat{\h}_{2,2}^{\mathrm{MC}}}
\newcommand{\hzerotwoMC}{\widehat{\h}_{0,2}^{\mathrm{MC}}}
\newcommand{\htwozeroMC}{\widehat{\h}_{2,0}^{\mathrm{MC}}}
\newcommand{\vloneone}{\widehat{\mathbb{V}}_{l,1,1}^{\mathrm{MC}}}
\newcommand{\voneoneMC}{\widehat{\mathbb{V}}_{1,1}^{\mathrm{MC}}}
\newcommand{\p}{p}
\newcommand{\q}{q}
\newcommand{\rr}{r}
\newcommand{\vech}{\mathrm{vech}}
\newcommand{\feq}[1]{Eq.~(\ref{#1})} 
\newcommand{\feqs}[2]{Eqs.~(\ref{#1}) and (\ref{#2})} 
\newcommand{\feqss}[3]{Eqs.~(\ref{#1}), (\ref{#2}) and (\ref{#3})} 
\newcommand{\fsec}[1]{Section~\ref{#1}}
\newcommand{\fsecs}[2]{Sections~\ref{#1} and \ref{#2}} 

\newcommand{\ffig}[1]{Fig.~\ref{#1}} 
\newcommand{\ftbl}[1]{Table~\ref{#1}}

\newcommand{\BRrefeq}[1]{Eq.~(\ref{#1})}

\newcommand{\prop}[1]{Proposition.~\ref{#1}} 

\title{Covariance estimation using h-statistics in Monte Carlo and multilevel Monte Carlo methods}
\author{Sharana Kumar Shivanand\textsuperscript{a,b}\thanks{{E-mail:} \texttt{sshivanand@turing.ac.uk}} \\
\small \textsuperscript{a} The Alan Turing Institute, London, United Kingdom \\
\small \textsuperscript{b} Department of Engineering, University of Cambridge, United Kingdom}
\date{\today}

\newenvironment{keyword}{%
	\noindent \emph{Keywords:}}{\newline}

\begin{document}
\sloppy
\maketitle


\begin{abstract}
We present novel Monte Carlo (MC) and multilevel Monte Carlo (MLMC) methods to determine the unbiased covariance of random variables using h-statistics. The advantage of this procedure lies in the unbiased construction of the estimator's mean square error in a closed form. This is in contrast to conventional MC and MLMC covariance estimators, which are based on biased mean square errors defined solely by upper bounds, particularly within the MLMC. The numerical results of the algorithms are demonstrated by estimating the covariance of the stochastic response of a simple 1D stochastic elliptic PDE such as Poisson's model.
\end{abstract}


\begin{keyword}
	\textit{covariance estimator, Monte Carlo, multilevel Monte Carlo, h-statistics, uncertainty quantification}
\end{keyword}

\section{Introduction}

\rvsi{Covariance estimation of random variables plays a fundamental role in a wide range of real-world applications, ranging from finance for portfolio optimization \cite{ledoit_power_2022} to engineering for system control \cite{hotz_covariance_1987} and astrophysics for orbit propagation \cite{cano_covariance_2023}. 
Furthermore, it is well known that determination of covariance is a fundamental task in uncertainty quantification (UQ). In forward UQ problems, it is crucial to approximate the Gaussian fields or processes, for instance, in the context of second-moment analysis techniques \cite{liu_random_1986,haldar2000probability}. 
This importance becomes apparent in Bayesian inverse problems, where the covariance derived from the forward problem plays a key role in parameter estimation through Kalman filtering-type data assimilation \cite{rosic_parameter_2013,rosic_sampling-free_2012,matthies_inverse_2016}. 
It is also instrumental in tasks related to uncertainty analysis, encompassing both variance- \cite{sobol1993sensitivity,sobol_global_2001} and covariance-based \cite{liu_generalized_2021,mycek_multilevel_2019} sensitivity analysis.} 
Additionally, covariance estimation has significant relevance in data science and machine-learning applications, including principal component analysis for dimensionality reduction \cite{jolliffe_principal_2016}, metric learning \cite{zadeh_geometric_2016,maurais_multi-fidelity_2023}, and spatial statistics \cite{cressie2015statistics}.

In this study, we focus on obtaining the covariance estimate of random data generated by solving physics-based stochastic differential equations, with an emphasis on forward UQ. This task involves propagating the input randomness to the output stochastic response.  
Various algorithms have been developed to address this challenge, typically categorized into functional approximation methods and the direct integration of statistics, as extensively discussed in \cite{xiu_numerical_2010, sullivan_introduction_2015, le_maitre_spectral_2010}.

While functional approximation methods are sensitive to the magnitude of input dimensionality, this study concentrates on direct integration techniques, particularly random sampling-based approaches. Among these techniques, the Monte Carlo (MC) method has traditionally been recognized as the gold standard for solving stochastic problems because of its simplicity and resistance to the curse of dimensionality \cite{metropolis_monte_1949, fishman_monte_1996, graham_stochastic_2013}. However, the MC can be computationally demanding and slow to converge.
To address these limitations, we explore the sampling-based Multilevel Monte Carlo method (MLMC) as an efficient alternative to MC.
Under suitable conditions, MLMC reduces computational costs by hierarchically adjusting the number of stochastic samples as model fidelity increases, making it a compelling option for practical applications. 

To the best of authors' knowledge, the concept of MLMC was introduced by \cite{Heinrich2001} to estimate high-dimensional parameter-dependent integrals. It gained prominence and wide application in computational finance, particularly in solving It\^{o}'s stochastic ordinary differential equations, as popularized by \cite{giles2008}. However, earlier studies primarily focused on approximating the sample mean of the output response and lacked a comprehensive characterization of the probabilistic solution. Subsequently, \cite{mishra_multi-level_2012} developed an unbiased MLMC sample variance estimator, which was further analyzed by \cite{bierig_convergence_2015} in the context of elliptic random obstacle problems. An alternative MLMC variance estimator based on h-statistics was introduced by \cite{krumscheid_quantifying_2020}. In addition, the estimation of higher-order moments, such as skewness and kurtosis, using MLMC has been explored in \cite{bierig_estimation_2016} and \cite{krumscheid_quantifying_2020}.

Early literature on MLMC covariance estimation primarily focused on data assimilation to address Bayesian inverse problems. For example, to the best of our knowledge,  \cite{hoel_multilevel_2016} pioneered the development of MLMC covariance estimation using filtering methods to address stochastic time-dependent problems. Building upon this, \cite{chernov_multilevel_2021}  extended the approach to spatio-temporal systems.
In addition, \cite{popov_multifidelity_2021} introduced an extension related to the multifidelity ensemble Kalman filter algorithm, leveraging the theory of multivariate control variates and reduced order modelling. 
In \cite{mycek_multilevel_2019}, the MLMC theory was further extended, particularly for covariance estimation, and the resulting multilevel estimator was applied to sensitivity analysis to derive multilevel estimates of Sobol\textquoteright\ indices.

\rvso{As to the challenge of handling large covariance matrices in a MLMC context, \cite{dolz_data_2023} propose data sparse multilevel sample covariance estimator, which use $\mathcal{H}^2$-approximations to compress sample covariances on each level. Also, the authors in \cite{chernov_simple_2023} approximate the covariance function in a multilevel fashion, and additionally derive sparse tensor product approximation variants to overcome the curse of dimensionality.}

More recently, the issue of covariance matrix definiteness within generalized variance-based multifidelity frameworks was addressed in \cite{maurais_multi-fidelity_2023, maurais_multifidelity_reg_2023}. 
In \cite{maurais_multi-fidelity_2023}, an estimator based on the principles of log-Euclidean geometry for symmetric positive definite (SPD) matrices was introduced.
On the other hand, \cite{maurais_multifidelity_reg_2023} outlined a multifidelity covariance estimator as the solution to a regression problem formulated on the tangent spaces of the product manifolds of SPD matrices. 
Both methodologies guarantee matrix definiteness through their construction.
However, in the case of the log-Euclidean multifidelity estimator, as presented in \cite{maurais_multifidelity_reg_2023}, the derivation of the mean square error (MSE) remains a challenging task, and the available solution is presently limited to the first order in the log-Euclidean metric.
By contrast, the regression estimator in \cite{maurais_multifidelity_reg_2023} does not provide a performance guarantee. This implies that no error bounds that delineate the dependence on computational cost have been elucidated. Moreover, there is no assurance that the resulting regression covariance estimator is bias-free.
It is imperative to emphasize that our study, in contrast to the aforementioned research, does not encompass broader multifidelity frameworks concerning covariance matrices. Instead, it fundamentally centers on a Monte Carlo-based scalar covariance estimation.

Traditional MC covariance estimation relies on a closed-form yet biased sampling error \cite{bierig_convergence_2015}, whereas MLMC covariance estimators depend on biased and worst-case-bounded error estimates, as discussed in \cite{hoel_multilevel_2016,mycek_multilevel_2019}. This research is dedicated to improving the overall accuracy and efficiency of MLMC covariance estimation. Our objective is to introduce a novel approach that employs h-statistics \cite{dwyer1937,halmos1946} to estimate the unbiased covariance of random variables, \rvsn{with minimal variance}, using both the MC and MLMC methods. 
\rvsn{The mean-square error estimates that result from this allow us to derive unbiased, closed-form versions that are potentially more accurate than the error bounds corresponding to conventional estimators.}
Although these calculations can be complex, we harness the powerful symbolic capabilities of the \textit{mathStatica} library \cite{Colin_Rose:2002} in \textit{Mathematica}. \rvsn{The equations related to h-statistics based MC and MLMC covariance estimation are generated in a Mathematica notebook \cite{shivanandGitHub2024}}. Notably, \cite{krumscheid_quantifying_2020} pioneered the successful implementation of h-statistics to establish unbiased and closed-form error estimates for MC and MLMC central moments, spanning from the second to fourth order.

We evaluate the effectiveness of our proposed methods by applying them to linear elliptic problems, as exemplified by Poisson's model. As a specific case, we consider a straightforward 1D steady-state heat equation, where uncertainty is incorporated into the material properties, particularly in the form of spatially constant random fields that affect the thermal conductivity. This study focuses exclusively on material uncertainties, assuming that all other model parameters are deterministic. However, it is important to note that our numerical approach is versatile and can be extended to handle various types of uncertainty. Both the MC and MLMC methods propagate this uncertainty, enabling the estimation of key statistics including the covariance, mean, and variance of the temperature field.

The remainder of this paper is organized as follows. In \fsec{UncertaintyQuantification}, we describe the problem. The theoretical procedure for both conventional and h-statistics-based MC covariance estimates is elaborated in \fsec{Sec:MC}. The corresponding MLMC covariance estimators are detailed in \fsec{Sec:MLMC}. The deterministic and stochastic settings of Poisson's model are provided in  \fsec{Sec:application}. In \fsec{NumericalResults}, the numerical results of the MLMC and MC estimators, when implemented in a one-dimensional steady-state heat conduction example, are explained. Finally, conclusions are drawn in \fsec{Conclusion}.

\section{Problem description}\label{UncertaintyQuantification}

Consider a physical system that occupies the spatial domain $\spatialdomain\subset\Euclideanspace$ in a d-dimensional Euclidean space and is represented by an abstract equilibrium equation,
\begin{equation}\label{abs_main_eq}
	\mathcal{A}(q(x),\uu(x))=f(x).
\end{equation}
For simplicity, $\uu(x) \in \U$ characterizes the system's state at spatial location $x\in\spatialdomain$ in the Hilbert space $\U$, $\mathcal{A}$ is a (potentially non-linear) operator describing the system's physics, and $f\in \U^*$ is an external influence (action, excitation, or loading). 
We further assume that the model depends on the parameter set $q \in \mathcal{Q}$ and is accompanied by a suitable boundary and/or initial conditions.
It should be noted that the model in the previous equation only describes the physical system in the spatial domain $x$ for conciseness, despite the fact that it can be generalized to time-dependent processes \rvsn{\cite{hoel_multilevel_2016, chernov_multilevel_2021}}.

In reality, complete knowledge of model parameters in \feq{abs_main_eq} is not attainable; in other words, uncertainty may arise as a result of randomness in external influence $f$,  initial or boundary conditions, geometry $\spatialdomain$, and the coefficient of operator $\mathcal{A},$ that is, parameter $q$. 
Although the theory described in the following sections is not dependent on this choice and is generic enough to accommodate all the specified (single or combination of) scenarios, this study concentrates on introducing randomness in coefficient $q$. For instance, in the theory of steady-state heat transfer, the quantity $q$ denotes thermal conductivity, as seen in \fsec{Sec:application}.
Here, we assume that $q$ is modelled as a spatially constant random field $q(\omega)$ with finite second-order moments on a probability space $(\varOmega, \mathfrak{F}, \mathbb{P})$, where $\varOmega$ denotes the set of elementary events, $\mathfrak{F}$ is a $\sigma$-algebra of events, and $\mathbb{P}$ is a probability measure.
Following this, \BRrefeq{abs_main_eq} is rewritten in a stochastic form:
\begin{equation}\label{basic_eq}
	\mathcal{A}(q(\omega),\uu(x,\omega))=f(x),
\end{equation}
which must then be solved for $\uu(x,\omega) \in \U\otimes L_2(\varOmega, \mathfrak{F}, \mathbb{P})$. 

\rvso{Because of the spatial and stochastic dependence of the solution $\uu \equiv\uu(x,\omega)$ in the previous equation, it is first discretized in a spatial domain, that is,~we search for the solution in a finite subspace $\U^h\subset \U$.  After rewriting the problem in a variational form, the spatial discretization $\uh(x, \omega)\in \U^h$---$h$ is the discretization parameter---can take, for example, the finite element form  on a sufficiently fine spatial grid\footnote{Other types of discretization techniques can be considered as well}.}

\rvso{The goal of this study is to determine the covariance between the finite variance random variables $u^h_i \equiv \uh(x_i, \event)$ and $u^h_j \equiv \uh(x_j, \event)$ at spatial points $x_i$ and $x_j$, respectively i.e.,
\begin{equation}\label{Eq:cov_analytical}
	\cov(u^h_i, u^h_j) = \mathbb{E}((u^h_i - \mathbb{E}(u^h_i)) (u^h_j - \mathbb{E}(u^h_j))).
\end{equation}
Here, the expectation operator $\mathbb{E}$ or the mean $\mu$ of a scalar-valued function, let's say of a random field $Y(x, \event): \spatialdomain \times \samplespace \rightarrow \mathbb{R} $, is defined analytically in the form
\begin{equation}\label{math_Exp_Y}
	\mathbb{E}(Y) = \mu(Y)  = \int_{\samplespace} Y \, \mathbb{P}(\textrm d\event),
\end{equation}
where $Y$ denotes $u^h_i$ or $u^h_j$ in \feq{Eq:cov_analytical}.}

Note that $\cov(u^h_i, u^h_j)$ is the $(i, j)$ entry of the covariance matrix
\begin{equation}\label{Eq:cov_analytical_vec}
	\cov(\vek{U}^h) = \mathbb{E}((\vek{U}^h - \mathbb{E}(\vek{U}^h)) (\vek{U}^h - \mathbb{E}(\vek{U}^h))^{\mathrm{T}}),
\end{equation}
where $\vek{U}^h(\event) \equiv \vek{U}^h := [u^h_1, u^h_2,....,u^h_M]^{\mathrm{T}}  \in \mathbb{R}^\nodes$ denotes a $\nodes$-dimensional random solution vector.

 Furthermore, $\cov(\vek{U}^h)$ takes values in the open convex cone of real-valued second-order positive-definite tensor,
\begin{equation}\label{Eq:PosSymSpace}
	\PosSymMatSpace(\nodes):=\{\vek{H}\in (\mathbb{R}^\nodes\otimes\mathbb{R}^\nodes) \mid 
	\vek{H}=\vek{H}^T, \,\vek{z}^T \vek{H}  \vek{z}>\vek{0}, \forall 
	\vek{z}\in\mathbb{R}^\nodes\setminus\vek{0} \}.
\end{equation}

\section{Monte Carlo method}\label{Sec:MC}
As the analytical integration of \feq{math_Exp_Y} becomes a formidable challenge, numerical integration methods must be embraced. In this regard, the Monte Carlo (MC)  \cite{metropolis_monte_1949,fishman_monte_1996} method has traditionally emerged as the most common choice.
\subsection{Classical Monte Carlo estimator}\label{Sec:classMC}
To estimate the mean using the Monte Carlo procedure (MC), denoted as
\begin{equation}
	\mathbb{E}(\uh)  \approx \muMC(\uh) = \frac{1}{N} \sum_{k=1}^{N} \uh(x, \event_k),
\end{equation}
we rely on a set of $\samplesize>1$ independent sample points $\left\lbrace \event_1, \event_2,...,\event_N \right\rbrace \in \samplespace $ drawn from the random field $\uh(x, \event)$. Consequently, the computation of the covariance in \feq{Eq:cov_analytical} is transformed into the following expression:
\begin{equation}\label{Eq:cov_sample}
	\covMC\uhij= \frac{N}{N-1} \muMC((u^h_i - \muMC(u^h_i)) (u^h_j - \muMC(u^h_j))),
\end{equation}
where the expectation operator $\mathbb{E}$ is replaced by $\muMC$. Note that in the previous equation, $\covMC\uhij$ serves as an unbiased estimator, implying that $\mathbb{E}(\covMC\uhij) = \cov\uhij$, owing to the presence of Bessel's correction term $N/N-1$.

Subsequently, we express the mean squared error (MSE) of the MC covariance estimator as
\begin{equation}\label{Eq:MSEcov}
	\textrm{MSE}(\covMC\uhij) = \var(\covMC\uhij) + (\cov\uhij- \cov\uij)^2.
\end{equation}
The first term reflects sampling or stochastic accuracy. The variance is defined as
\begin{equation}\label{Eq:varcovMC}
	\var(Z) = \mathbb{E}((Z-\mathbb{E}(Z))^2) = \int_\samplespace (Z-\mathbb{E}(Z))^2 \probability(\text{d}\event),
\end{equation}
where $Z=\covMC\uhij$.
On the other hand, the second term in \feq{Eq:MSEcov} represents the square of the deterministic discretization error, primarily in the spatial domain. Therefore, to attain an overall mean squared error of $\epsilon^2$, one may ensure that both terms are less than $\epsilon^2/2$. However, such an equal splitting of accuracy is not a requirement and can be set otherwise.

This study specifically focuses on analyzing the sampling error, which can be further expanded (as shown in \cite{bierig_convergence_2015,mycek_multilevel_2019}) as follows:
\begin{equation}\label{Eq:varMCcov}
	\var(\covMC\uhij) = \frac{\mu_{2,2}\uhij}{\N} - \frac{(\N-2)\cov\uhij^2}{\N(\N-1)}+\frac{\var(u^h_i) \var(u^h_j)}{\N(\N-1)}. 
\end{equation}
Here $\mu_{2,2}\uhij$ denotes the bivariate central moment:
\begin{equation}\label{Eq:bivarrs}
	\mu_{\p,\q}\uhij=  \mathbb{E}((u^h_i - \mathbb{E}(u^h_i))^\p (u^h_j - \mathbb{E}(u^h_j))^\q),
\end{equation} with orders $\p,\q$ belonging to a set of non-negative integers $\mathbb{Z}_{\geq 0}$.

It is noteworthy that the variance, as outlined in \feq{Eq:varMCcov}, remains analytically tractable. Consequently, to attain better accuracy in this estimation, the challenge lies in obtaining an unbiased estimator, potentially in a closed form.
Therefore, this study formulates an unbiased framework for the approximation of sampling errors in the MC covariance estimator using h-statistics \cite{dwyer1937}.


\subsection{Monte Carlo estimation by h-statistics}

\subsubsection{h-statistics: unbiased estimation of central moments}
The $\p^{\text{th}}$ univariate central moment of solution $\uh$ can be expressed analytically as follows:
\begin{equation}\label{math_Exp_h}
	\mu_{\p}(\uh) =	\mathbb{E}((\uh-\mathbb{E}(\uh))^{\p}) = \int_\samplespace (\uh-\mathbb{E}(\uh))^\p \probability(\text{d}\event).
\end{equation}
Here, $\mathbb{E}(\uh)$ represents the mean or first raw moment (not to be confused with the first central moment, which equals zero when $\p=1$). When $\p=2$, $\var(\uh):=\mathbb{E}\left((\uh-\mathbb{E}(\uh))^2\right)$ corresponds to the second central moment or variance as described in \feq{Eq:varcovMC}.

The h-statistics, denoted by $\hrMC \equiv \hrMC(\uh)$, serves as an unbiased estimator of the central population moments $\mu_\p(\uh)$ in \feq{math_Exp_h} \cite{dwyer1937, Colin_Rose:2002}. The \textquoteleft MC\textquoteright\ superscript of $\widehat{\h}_\p$ signifies that these statistics are obtained through random Monte Carlo-based sampling.  For example, the first two h-statistics can be defined as
\begin{equation}\label{Eq:h-stat}
	\begin{aligned}
		\mu_1 \approx \widehat{\h}_1^{\text{MC}}  &= 0, \\
		\mu_2 \approx \htwoMC &= \frac{N\s_2-\s_1^2}{N(N-1)}, \\
	\end{aligned}
\end{equation}
in which $\s_1$ and $\s_2$ are the power sums, given as 
\begin{equation}
	\s_a(\uh) := \sum_{k=1}^{\N}\uh(x, \event_k)^a,
\end{equation} for $a\in\mathbb{Z}_{\geq 0}$. Using the power series representation of h-statistics, the authors in \cite{Colin_Rose:2002} developed a \textit{Mathematica}-based package called \textit{mathStatica}, which efficiently generates h-statistics for any value of $\p$.

What distinguishes h-statistics from other unbiased estimators are its properties, such as \cite{halmos1946}: 
\begin{enumerate}
	\item \textbf{Unbiasedness}: The expectation of the h-statistic equals the corresponding population central moment, i.e., $\mathbb{E}(\hrMC) = \mu_\p$.
	\item \textbf{Symmetry}: Among all unbiased estimators of $\mu_\p$, $\hrMC$ is the only one that exhibits symmetry\footnote{A function or an estimate is considered symmetric when it is not influenced by the order in which observations are considered.}.
	\item \textbf{Minimum variance}: Out of all unbiased estimators of $\mu_\p$, $\hrMC$ stands out for its minimal variance, denoted as $\var(\hrMC)$.
\end{enumerate}


\subsubsection{Monte Carlo estimation of covariance by h-statistics}\label{Sec:MChstat}

Covariance can be defined as a bivariate central moment with orders \(p, q = 1\). Thus, \feq{Eq:bivarrs} can be expressed as 
\begin{equation}
	\mu_{1,1}\uhij=  \mathbb{E}((u^h_i - \mathbb{E}(u^h_i)) (u^h_j - \mathbb{E}(u^h_j))). 
\end{equation}
This equation is identical to that of \feq{Eq:cov_analytical}, \rvsi{and the corresponding
bivariate h-statistics estimation of covariance is  \cite{Colin_Rose:2002}}
\begin{equation}\label{Eq:h11}
	\mu_{1,1}\uhij\approx \honeoneMC\uhij=	\frac{\N \s_{1,1}-\s_{0,1} \s_{1,0}}{(\N-1) \N}.
\end{equation}
Here, \(\s_{1,1}\), \(\s_{0,1}\), and \(\s_{1,0}\) represent the bivariate power sum,
\begin{equation}
	\s_{a,b}\uhij= \sum_{k=1}^{\N} u^h_i(x, \event_k)^a u^h_j(x, \event_k)^b,
\end{equation} 
where \(a, b \in \mathbb{Z}_{\geq 0}\).

As a result, the MSE of $\honeoneMC\equiv \honeoneMC\uhij$ can be rewritten from \feq{Eq:MSEcov} as,
\begin{equation}\label{Eq:MSEhstat}
	\textrm{MSE}(\honeoneMC) = \var(\honeoneMC) + (\mu_{1,1}\uhij- \mu_{1,1}\uij)^2,
\end{equation}
where the sampling accuracy is recast as bivariate population central moments:
\begin{equation}\label{Eq:varMCcovmu}
	\var(\honeoneMC) = \frac{\mu_{2,2}\uhij}{\N} - \frac{(\N-2)\mu_{1,1}\uhij^2}{\N(\N-1)}+\frac{\mu_{0,2}\uhij\mu_{2,0}\uhij}{\N(\N-1)}. 
\end{equation}
Here, $\mu_{0,2}\uhij= \var(u^h_j)$ denotes the variance of random variable $u^h_j$, and $\mu_{2,0}\uhij= \var(u^h_i)$ is the variance of $u^h_i$.

It is important to note that replacing the bivariate central moments with h-statistics in the previous equation does not directly yield an unbiased estimate because of the inherent nonlinear transformation of statistics. For simplicity and clarity, we drop the random variable pair notation $\uhij$ with h-statistics (and the population central moments hereafter). Consequently, following the approach in \cite{krumscheid_quantifying_2020}, we first consider the Monte Carlo estimate of $\var(\honeoneMC)$ in ansatz form:
\begin{equation}\label{Eq:varMCcovh21}
	\var(\honeoneMC) \approx \varMC(\honeoneMC) = {\aone \htwotwoMC} + {\atwo (\honeoneMC)^2 + \athree \hzerotwoMC \htwozeroMC},
\end{equation}
where $\aone, \atwo$ and $\athree$ are the scalar and real-valued coefficients. To obtain an unbiased estimator $\varMC(\honeoneMC)$, the condition $\mathbb{E}(\varMC(\honeoneMC)) = \var(\honeoneMC)$ must hold. 
\rvsi{Furthermore, one can
determine the population moment of a sample moment and represent it in either the central or raw moment, as described in \cite{Colin_Rose:2002}. Following this, we} calculate the expectation of $\varMC(\honeoneMC)$ in terms of central moments as
\begin{align}\label{Eq:varMCcovh22}
	\begin{split}
		\mathbb{E}(\varMC(\honeoneMC)) =& \frac{(\aone\N+\atwo+\athree) \mu_{2,2}}{\N} + \frac{(2 \athree + \atwo (2 - 2 \N + 				
			\N^2)){\mu_{1,1}}^2}{\N(\N-1)}+ \\ 
		&	\frac{(\atwo + \athree (-1 + \N)^2)\mu_{0,2} \mu_{2,0}}{\N(\N-1)}.
	\end{split}
\end{align}
By comparing the above equation with \feq{Eq:varMCcovmu} and solving a symbolic \rvsi{linear system of equations}, the values of the constants can be determined to be $\aone = \frac{\N-1}{3-2\N+\N^2}$, $\atwo = \frac{-5+4\N-\N^2}{-6+7\N-4\N^2+\N^3}$ and $\athree=\frac{\N-1}{-6+7\N-4\N^2+\N^3}$.

Finally, by substituting the coefficients into \feq{Eq:varMCcovh21} and expanding the h-statistics in the power sums, we obtain the unbiased closed-form expression of $\varMC(\honeoneMC)$ as follows:
\begin{align}\label{Eq:varMCh11}
	\begin{split}
		\varMC(\honeoneMC)  = & \frac{1}{(\N-3) (\N-2) (\N-1)^2 \N^2} \left( \N \left(\left(-\N^2+\N+2\right) \s_{1,1}^2+ \right. \right. \\
		& \left. \left. (\N-1)^2 \left(\N \s_{2,2}-2 \s_{1,0} \s_{1,2}\right)+ (\N-1) \s_{0,2} \left(\s_{1,0}^2-\s_{2,0}\right)\right)+ \right. \\ 
		& \left.  \s_{0,1}^2 \left((6-4 \N) \s_{1,0}^2+(\N-1) \N \s_{2,0}\right)  \right. \\
		&  \left.   -2 \N \s_{0,1} \left((\N-1)^2 \s_{2,1}+(5-3 \N) \s_{1,0} \s_{1,1}\right)\right).
	\end{split}
\end{align}
\textcolor{black}{Note that, for the statistic $\varMC(\honeoneMC)$ to be well-defined, it is important to ensure a minimum sample size of \rvsi{$\samplesize>3$} is considered.}

In the interest of brevity, we express the previous equation as \(\varMC(\honeoneMC) = \voneoneMC/\N\), which is equivalent to
\begin{equation}\label{Eq:VoneoneMC}
	\voneoneMC = \N\, \varMC(\honeoneMC).
\end{equation} 
This implies that the stochastic convergence of the MC covariance estimator follows the order \(\mathcal{O}(\N^{-1})\). Consequently, the reduction in MC error occurs at a relatively slow rate, making it a computationally expensive method. Furthermore, the procedure becomes infeasible when implemented in systems and processes represented by high-fidelity models, which are commonly used in many science and engineering applications. For instance, decreasing the value of the semi-discretization parameter \(h\) increases the deterministic fidelity in \feq{basic_eq}, but also significantly increases the computational cost required to determine a single sample of the solution \(\uh(x,\event)\).

Following this, we outline the total computational cost of the MC covariance estimation for a given MSE of \(\epsilon^2\) in the following proposition:
\begin{proposition}
	Let us consider $c_\alpha$, $c_\gamma$, $\alpha$, and $\gamma$ to be positive constants. The error bounds are defined as follows:
	\begin{equation*}\label{eq:errorvarmc}
		\begin{split}
			(i)&\quad \text{The deterministic error decays as } |\mu_{1,1}\uhij - \mu_{1,1}\uij| \le c_\alpha h^{\alpha},\\
			(ii)&\quad \text{The computational cost to determine one sample of $\uh(x, \event)$ is given by} \\
			&\qquad\qquad \mathcal{C}(\uh) \le c_\gamma h^{-\gamma}.
		\end{split}
	\end{equation*}
	Then, for any $0 < \epsilon < e^{-1}$, the Monte Carlo (MC) covariance estimator $\honeoneMC$ with $\samplesize = \mathcal{O}(\epsilon^{-2})$ and $\elesize = \mathcal{O}(\epsilon^{1/\alpha})$ satisfies $\mathrm{MSE}(\honeoneMC) < \epsilon^2$. Therefore, the corresponding MC computational cost to estimate the covariance is bounded by
	\begin{equation*}
		\mathcal{C}(\honeoneMC) \le c \epsilon^{-2 - \gamma/\alpha},
	\end{equation*}
	where $c > 0$.
\end{proposition}

\section{Multilevel Monte Carlo method}\label{Sec:MLMC}
To address the limitations of the MC method, we focus on the sampling-based multilevel Monte Carlo method (MLMC)  \cite{Heinrich2001,giles2008,giles_2015}. This approach offers a promising solution by distributing the sampling strategy across a hierarchy of fidelity levels, primarily in a sequential increase in fidelity from low to high. Under appropriate conditions (to be discussed in \fsec{Sec:MLcovhstat}), the resulting reduction in the number of stochastic samples as fidelity increases marks a significant increase in overall computational efficiency.

\subsection{Classical multilevel Monte Carlo estimator}
Consider a sequence denoted as $\left\lbrace \meshlevel=0,1,2...,\LL \right\rbrace$, which represents an increasing hierarchy of nested mesh grids $\mathcal{P}_\meshlevel$ within the context of a decreasing element size $h$. These grids form a regular nondegenerate partition of the computational domain $\mathcal{G}$, as defined in \feq{abs_main_eq}. Here, $\meshlevel$ corresponds to the mesh level and $\LL$ represents the finest available mesh level.

The primary objective of the Multilevel Monte Carlo method (\acrshort{mlmc}) is to ascertain statistical properties, such as covariance \cite{hoel_multilevel_2016,mycek_multilevel_2019}, of the solution $\uhL(\event)$ on the finest mesh level $\LL$:
\begin{equation}\label{Eq:covMLMC}
	\covML(u^{h_L}_i, u^{h_L}_j)_{\left\lbrace \N_l\right\rbrace} := \covMC(u^{h_0}_i, u^{h_0}_j)_{\N_0} + \sum_{\meshlevel=1}^{\numofmesh}  \covMC(u^{h_l}_i, u^{h_l}_j)_{\N_l} - \covMC(u^{h_{l-1}}_i, u^{h_{l-1}}_j)_{\N_l}.
\end{equation}
Here, $\Nl$ denotes the number of samples at each mesh level and $\left\lbrace \Nl \right\rbrace:=\left\lbrace \N_0,\N_1,...,\N_L \right\rbrace$ represents the set of samples required to determine $\covML$. In addition, $ \covMC(u^{h_0}_i, u^{h_0}_j)_{\N_0}$ is the \acrshort{mc} estimator of $\cov(u^{h_0}_i, u^{h_0}_j)$ at level $l=0$ using $N_0$ samples. Similarly, $\covMC(u^{h_l}_i, u^{h_l}_j)_{\N_l}$ and $\covMC(u^{h_{l-1}}_i, u^{h_{l-1}}_j)_{\N_l}$ represent the approximation of $\cov(u^{h_l}_i, u^{h_l}_j)$ and $\cov(u^{h_{l-1}}_i, u^{h_{l-1}}_j)$ with $\meshlevel>0$ and $\Nl$ samples, respectively. 
To simplify the discussion, we introduce the following:
\begin{equation}\label{Eq:Zl}
	\Yl= 
	\begin{cases}
		\covMC(u^{h_0}_i, u^{h_0}_j)_{\N_0},\quad &\meshlevel=0,\\
		\covMC(u^{h_l}_i, u^{h_l}_j)_{\N_l} - \covMC(u^{h_{l-1}}_i, u^{h_{l-1}}_j)_{\N_l},\quad &\meshlevel>0.
	\end{cases}
\end{equation}
Following this, the covariance estimator $\covML$ from \feq{Eq:covMLMC} can be expressed as
\begin{equation}\label{Eq:sumZl}
	\covML	= \sum_{\meshlevel=0}^{\numofmesh} \Yl.
\end{equation}
It is important to note that each term $\Yl, \ \meshlevel\ge{0}$ is sampled independently, and hence remains uncorrelated. However, for $l>0$, the quantities $(u^{h_l}_i, u^{h_l}_j)_{\N_l}$ and $(u^{h_{l-1}}_i, u^{h_{l-1}}_j)_{\N_l}$ in $\Yl$ are strongly correlated and originate from the same random seeds.

A multilevel covariance estimator, represented by \feq{Eq:covMLMC} or, in alternative terms, the sample multilevel covariance matrix $\covML(\vek{U}^{h_L})$ of the form shown in \feq{Eq:cov_analytical_vec}, has a limitation in that it does not guarantee the positive-definiteness property required for covariance matrices, as explained in \feq{Eq:PosSymSpace}.
This is because subtracting one positive-definite matrix from another does not guarantee that the resulting matrix is positive-definite, which is vital to the stability of the algorithm \cite{hoel_multilevel_2016,maurais_multi-fidelity_2023}.
To overcome this drawback, we follow a simple post-hoc procedure (as also depicted in \cite{hoel_multilevel_2016}) by performing eigen-decomposition on
\begin{equation}
	\covML = \sum_{i=1}^{\nodes} \lambda_i \vek{q}_i \vek{q}_i^T,
\end{equation}
where $\lambda_i\in\mathbb{R}$ denotes the eigenvalues, and $\vek{q}_i\in\mathbb{R}^{\nodes}$ are the corresponding eigenvectors. 
Positive definiteness is ensured by removing the zero- and negative-valued eigenvalues ($\lambda_i\le0$), that is
\begin{equation}
	\mrm{c}\widetilde{\mrm{o}}\mrm{v}^{\textrm{ML}}  = \sum_{i=1;\lambda_i>0}^{\nodes} \lambda_i \vek{q}_i \vek{q}_i^T.
\end{equation}
As alternatives to spectral decomposition, which may involve less intrusiveness and/or reduced computational costs, various numerical techniques, including banding, shrinkage, thresholding, and localization (as outlined in \cite{hoel_multilevel_2016}), offer the assurance of positive-definiteness of the covariance matrix, particularly in higher dimensions. In contrast, \cite{maurais_multi-fidelity_2023} introduced a multifidelity covariance estimator that enforces the definiteness constraint by construction, utilizing the log-Euclidean mapping of symmetric and positive definite (SPD) matrices.

As the multilevel Monte Carlo covariance estimate $\covML$ on the finest grid $L$ is obtained as the telescopic sum of the difference of the \acrshort{mc} covariance estimates on the coarser grids, the MSE corresponding to \feq{Eq:sumZl} takes the following form:
\begin{equation}\label{Eq:MSEcovML}
	\textrm{MSE}(\covML)  =  \sum_{\meshlevel=0}^{\numofmesh} \var(\Yl) + (\cov(u^{h_L}_i, u^{h_L}_j) - \cov(u_i, u_j))^2.
\end{equation}
Similar to \feq{Eq:MSEcov}, the error estimate above is composed of two terms. The first denotes the variance of  the covariance estimator $\covML$; that is, $\var(\covML) = \var(\sum_{\meshlevel=0}^{\numofmesh}\Yl) = \sum_{\meshlevel=0}^{\numofmesh} \var(\Yl)$. On the other hand, the second part represents the square of the spatial discretization error. 
As emphasized by the authors in \cite{mycek_multilevel_2019}, it is essential to determine the upper bound of the sampling error to ensure the convergence of the multilevel estimator. Therefore, it is derived as
\begin{align}\label{eq:samperrul}
	\begin{split}
		\sum_{\meshlevel=0}^{\numofmesh} \var(\Yl) &\le \frac{1}{2} \sum_{\meshlevel=0}^{\numofmesh} \frac{1}{\Nl-1} \left( 
		\sqrt{\mu_{4}(X^{^-}(u^{h_l}_i))\mu_{4}(X^{^+}(u^{h_l}_j))} + \right. \\
		& \left. \sqrt{\mu_{4}(X^{^+}(u^{h_l}_i))\mu_{4}(X^{^-}(u^{h_l}_j))} \right),
	\end{split}
\end{align}
where $\mu_{4}$ denotes the fourth-order central moment (see \feq{math_Exp_h} for $r=4$). Furthermore, $X^{^+}(u^{h_l}_i) = u^{h_l}_i + u^{h_{l-1}}_i$  and $X^{^-}(u^{h_l}_i) = u^{h_l}_i - u^{h_{l-1}}_i$. 

\rvsn{As the previous equation is analytical, if one estimates the sample fourth-order moment as
$$ \mu_4 \approx \widehat{\mu}_4^{\mathrm{MC}} = \frac{1}{\Nl} \sum_{k=1}^{\Nl} (X_k - \muMC(X))^4, $$
where $X$ is $X^{^-}$ or $  X^{^+}$, we obtain the approximation for the sampling error. Meaning that,
\begin{equation}\label{eq:appsamperrul}
	\sum_{\meshlevel=0}^{\numofmesh} \var(\Yl) \approx \sum_{\meshlevel=0}^{\numofmesh} \varMC(\Yl),
\end{equation}
which remains biased and worst-case bounded.}

\subsection{Multilevel Monte Carlo estimator using h-statistics}\label{Sec:MLcovhstat}
As discussed in \fsec{Sec:MC}, the objective of this study is to obtain an unbiased estimate of the sampling error of the estimator $\covML$, which is also expressed in a closed form. To this, we rewrite the classical multilevel covariance estimator from \feq{Eq:covMLMC} using h-statistics as 
\begin{equation}\label{Eq:covMLMChstat}
	\honeoneML(u^{h_L}_i, u^{h_L}_j)_{\left\lbrace \N_l\right\rbrace} :=\honeoneMC(u^{h_0}_i, u^{h_0}_j)_{\N_0} + \sum_{\meshlevel=1}^{\numofmesh}  \honeoneMC(u^{h_l}_i, u^{h_l}_j)_{\N_l} - \honeoneMC(u^{h_{l-1}}_i, u^{h_{l-1}}_j)_{\N_l}.
\end{equation}
Because the second part of the above equation (i.e., for $\meshlevel>0$) is composed of variables $u^{h_l}_i, u^{h_l}_j, u^{h_{l-1}}_i$ and $ u^{h_{l-1}}_j$, one may equivalently re-represent it in the form
\begin{align}\label{Eq:covMLMChstat1}
	\begin{split}
		\honeoneML(u^{h_L}_i, u^{h_L}_j)_{\left\lbrace \N_l\right\rbrace} =&  \honeoneMC(u^{h_0}_i, u^{h_0}_j)_{\N_0} + \sum_{\meshlevel=1}^{\numofmesh}  \widehat{\h}_{1,1,0,0}^{\text{MC}}(u^{h_l}_i, u^{h_l}_j,u^{h_{l-1}}_i, u^{h_{l-1}}_j)_{\N_l} - \\
		& \widehat{\h}_{0,0,1,1}^{\text{MC}}(u^{h_l}_i, u^{h_l}_j,u^{h_{l-1}}_i, u^{h_{l-1}}_j)_{\N_l}.
	\end{split}
\end{align}
Here, the quadrivariate h-statistic $\widehat{\h}_{\p,\q, \rr,\s}^{\text{MC}}$ is an unbiased estimator of the corresponding central moment ${\mu}_{\p,\q, \rr,\s}$ with orders $\p,\q, \rr,\s \in\mathbb{Z}_{\geq 0}$ i.e.,
\begin{equation}
	\begin{split}
		\mu_{\p,\q, \rr,\s}(u^{h_l}_i, u^{h_l}_j, u^{h_{l-1}}_i, u^{h_{l-1}}_j) =  \mathbb{E}((u^{h_l}_i- \mathbb{E}(u^{h_l}_i))^\p (u^{h_l}_j - \mathbb{E}(u^{h_l}_j))^\q \\ (u^{h_{l-1}}_i - \mathbb{E}(u^{h_{l-1}}_i))^\rr (u^{h_{l-1}}_j - \mathbb{E}(u^{h_{l-1}}_j))^\s).
	\end{split}
\end{equation}
See \feq{Eq:bivarrs} for the definition of bivariate central moment.
Similar to \feq{Eq:Zl}, one may simplify by introducing
\begin{equation}\label{Eq:Zlhstat}
	\Zl= 
	\begin{cases}
		\honeoneMC(u^{h_0}_i, u^{h_0}_j)_{\N_0},\quad &\meshlevel=0,\\
		\widehat{\h}_{1,1,0,0}^{\text{MC}}(u^{h_l}_i, u^{h_l}_j,u^{h_{l-1}}_i, u^{h_{l-1}}_j)_{\N_l} - \widehat{\h}_{0,0,1,1}^{\text{MC}}(u^{h_l}_i, u^{h_l}_j,u^{h_{l-1}}_i, u^{h_{l-1}}_j)_{\N_l},\quad &\meshlevel>0,
	\end{cases}
\end{equation}
in \feq{Eq:covMLMChstat1}, which transforms to
\begin{equation}\label{Eq:hstatZl}
	\honeoneML	= \sum_{\meshlevel=0}^{\numofmesh} \Zl.
\end{equation}

For $\meshlevel=0$, $Z_0 = \honeoneMC$ is described by \feq{Eq:h11}. \rvsi{On the other hand, $\Zl$ for $\meshlevel>0$ is defined
as \cite{Colin_Rose:2002}}
\begin{equation}
	\Zl = \frac{\Nl s_{1,1,0,0}-s_{0,1,0,0} s_{1,0,0,0}}{(\Nl-1) \Nl}-\frac{\Nl s_{0,0,1,1}-s_{0,0,0,1} s_{0,0,1,0}}{(\Nl-1) \Nl}.
\end{equation}
Here, the quadrivariate power sum $\s_{a,b,c,d}$, for each $\ml$, is given as 
\begin{equation}
	\s_{a,b,c,d}(u^{h_l}_i, u^{h_l}_j,u^{h_{l-1}}_i, u^{h_{l-1}}_j) = \sum_{k=1}^{\N} u^{h_l}_i(\event_k)^a u^{h_l}_j(\event_k)^b u^{h_{l-1}}_i(\event_k)^c u^{h_{l-1}}_j(\event_k)^d,
\end{equation}
where $a, b, c, d$ are non-negative integers.

\subsubsection{MSE of h-statistics based MLMC covariance estimator}
Owing to the independent and uncorrelated sampling of $\Zl$ in \feq{Eq:hstatZl}, for each $\meshlevel$, the MSE of h-statistics based MLMC covariance is described as
\begin{equation}\label{Eq:MSEhstatML}
	\textrm{MSE}(\honeoneML)  = \sum_{\meshlevel=0}^{\numofmesh} \var(\Zl) + (\mu_{1,1}(u^{h_L}_i, u^{h_L}_j) - \mu_{1,1}(u_i, u_j))^2.
\end{equation}
Clearly, for $\meshlevel=0$, $\var(Z_0)$ is determined using \feq{Eq:varMCcovmu}. \rvsi{However, for $\ml>0$, 
the variance of the difference term $\Zl$ is written in the form of central moments:}
\begin{align}\label{Eq:varZl}
	\begin{split}
		\var(\Zl) =& -\frac{(\Nl-2) \mu _{0,0,1,1}^2}{(\Nl-1) \Nl}+\frac{2 \mu _{1,1,0,0} \mu _{0,0,1,1}}{\Nl}+\frac{\mu _{0,0,0,2} \mu _{0,0,2,0}}{(\Nl-1) \Nl}+\frac{\mu _{0,0,2,2}}{\Nl}+ \\
		& \frac{\mu _{0,2,0,0} \mu _{2,0,0,0}}{(\Nl-1) \Nl}+\frac{\mu _{2,2,0,0}}{\Nl}-\frac{2 \mu _{1,1,1,1}}{\Nl}-\frac{(\Nl-2) \mu _{1,1,0,0}^2}{(\Nl-1) \Nl}- \\
		& \frac{2 \mu _{0,1,1,0} \mu _{1,0,0,1}}{(\Nl-1) \Nl}-\frac{2 \mu _{0,1,0,1} \mu _{1,0,1,0}}{(\Nl-1) \Nl}.
	\end{split}
\end{align}
Evidently, the previous equation is analytical. Therefore, similar to the approach discussed in \fsec{Sec:MChstat}, to obtain a Monte Carlo based---unbiased---estimate of $\var(\Zl)$, an ansatz of the form
\begin{align}\label{Eq:varMCZl}
	\begin{split}
		\var(\Zl) \approx \varMC(\Zl) =& a_1 {({\widehat{\h}_{0,0,1,1}})^2}+a_2{ {\widehat{\h}_{1,1,0,0}} {\widehat{\h}_{0,0,1,1}}}+ a_3 {{\widehat{\h}_{0,0,0,2}} {\widehat{\h}_{0,0,2,0}}}+ a_4 {\widehat{\h}_{0,0,2,2}}+ \\
		& a_5 {{\widehat{\h}_{0,2,0,0}} {\widehat{\h}_{2,0,0,0}}}+ a_6 {\widehat{\h}_{2,2,0,0}}+ a_7 {{\widehat{\h}_{1,1,1,1}}}+ a_8 {({\widehat{\h}_{1,1,0,0}})^2}+ \\ 
		& a_9 {{\widehat{\h}_{0,1,1,0}} {\widehat{\h}_{1,0,0,1}}}+ a_{10}{{\widehat{\h}_{0,1,0,1}} {\widehat{\h}_{1,0,1,0}}},
	\end{split}
\end{align}
is assumed.
Note that here, we have dropped the repetitive superscript \textquoteleft $\text{MC}$\textquoteright\ of h-statistics for clarity. As the goal is to satisfy the condition, $\mathbb{E}(\varMC(\Zl)) = \var(\Zl)$, the expectation of $\varMC(\Zl)$ is determined in quadrivariate central moments as
\begin{align}
		\mathbb{E}(\varMC(\Zl)) =& \frac{\mu _{0,0,2,2} (a_1+a_2+a_3 \Nl)}{\Nl}+ \frac{\mu _{0,0,1,1}^2 \left(a_1 \Nl^2-2 a_1 \Nl+2 a_1+2 a_2\right)}{(\Nl-1) \Nl}+ \nonumber\\
		& \frac{\mu _{0,0,0,2} \mu _{0,0,2,0} \left(a_1+a_2 \Nl^2-2 a_2 \Nl+a_2\right)}{(\Nl-1) \Nl}+  \nonumber\\
		& \frac{\mu _{2,2,0,0} (a_{10} \Nl+a_7+a_9)}{\Nl}+\frac{\mu _{1,1,1,1} (a_4+a_5+a_6+a_8 \Nl)}{\Nl}+ \nonumber \\ 
		& \frac{\mu _{1,1,0,0} \mu _{0,0,1,1} \left(a_4+a_5+a_6 \Nl^2-2 a_6 \Nl+a_6\right)}{(\Nl-1) \Nl}+ \nonumber\\
		& \frac{\mu _{0,1,1,0} \mu _{1,0,0,1} \left(a_4 \Nl^2-2 a_4 \Nl+a_4+a_5+a_6\right)}{(\Nl-1) \Nl}+ \nonumber\\
		& \frac{\mu _{0,1,0,1} \mu _{1,0,1,0} \left(a_4+a_5 \Nl^2-2 a_5 \Nl+a_5+a_6\right)}{(\Nl-1) \Nl}+ \\
		& \frac{\mu _{1,1,0,0}^2 \left(a_7 \Nl^2-2 a_7 \Nl+2 a_7+2 a_9\right)}{(\Nl-1) \Nl}+ \nonumber\\
		& \frac{\mu _{0,2,0,0} \mu _{2,0,0,0} \left(a_7+a_9 \Nl^2-2 a_9 \Nl+a_9\right)}{(\Nl-1) \Nl}\nonumber.
\end{align}
\rvsi{By comparing the previous equation and \feq{Eq:varZl}, 
the subsequent linear system of equations is solved.} This results in scalar and real-valued constants as,  
$ a_1=-\frac{\Nl^2-4 \Nl+5}{\Nl^3-4 \Nl^2+7 \Nl-6}$, $a_2=-\frac{1-\Nl}{\Nl^3-4 \Nl^2+7 \Nl-6}$ $a_3=-\frac{1-\Nl}{\Nl^2-2 \Nl+3}$, $a_4=-\frac{2 (\Nl-1)}{\Nl^3-4 \Nl^2+7 \Nl-6}$, $a_5=-\frac{2 (\Nl-1)}{\Nl^3-4 \Nl^2+7 \Nl-6}$, $a_6=\frac{2 \left(\Nl^2-3 \Nl+4\right)}{\Nl^3-4 \Nl^2+7 \Nl-6}$, $a_7=-\frac{\Nl^2-4 \Nl+5}{\Nl^3-4 \Nl^2+7 \Nl-6}$, $a_8 =-\frac{2 (\Nl-1)}{\Nl^2-2 \Nl+3}$, $a_9 =-\frac{1-\Nl}{\Nl^3-4 \Nl^2+7 \Nl-6}$ and $a_{10}=-\frac{1-\Nl}{\Nl^2-2 \Nl+3}$.
Finally, by substituting these values in \feq{Eq:varMCZl} and expressing the corresponding h-statistics in the form of power sums, one finally gets the unbiased closed-form expression:
\begin{align}
		\varMC(\Zl)=& \frac{1}{(\Nl-3) (\Nl-2) (\Nl-1)^2 \Nl^2} \nonumber\\
		& \left( \s_{0,0,2,2} \Nl^4-2 \s_{1,1,1,1} \Nl^4+ \s_{2,2,0,0} \Nl^4-\s_{0,0,1,1}^2 \Nl^3-\s_{1,1,0,0}^2 \Nl^3- \right.\nonumber\\
		& \left. 2 \s_{0,0,1,0} \s_{0,0,1,2}\Nl^3- 2 \s_{0,0,2,2} \Nl^3+ 2 \s_{0,1,1,1} \s_{1,0,0,0} \Nl^3+2 \s_{0,1,0,0} \s_{1,0,1,1} \Nl^3+ \right.\nonumber\\
		& \left. 2 \s_{0,0,1,1} \s_{1,1,0,0} \Nl^3+ 2 \s_{0,0,1,0} \s_{1,1,0,1} \Nl^3+ 4 \s_{1,1,1,1} \Nl^3-2 \s_{1,0,0,0} \s_{1,2,0,0} \Nl^3-\right.\nonumber\\
		& \left. 2 \s_{0,1,0,0} \s_{2,1,0,0} \Nl^3-2 \s_{2,2,0,0} \Nl^3+  \s_{0,0,1,1}^2 \Nl^2+ \s_{0,2,0,0} \s_{1,0,0,0}^2 \Nl^2 + \right.\nonumber\\
		& \left. \s_{1,1,0,0}^2 \Nl^2+ 4 \s_{0,0,1,0} \s_{0,0,1,2} \Nl^2+ \s_{0,0,2,2} \Nl^2-4 \s_{0,0,1,1} \s_{0,1,0,0} \s_{1,0,0,0} \Nl^2-\right.\nonumber\\
		& \left. 2 \s_{0,0,1,0} \s_{0,1,0,1} \s_{1,0,0,0} \Nl^2- 4 \s_{0,1,1,1} \s_{1,0,0,0} \Nl^2- 2 \s_{0,0,1,0} \s_{0,1,0,0} \s_{1,0,0,1} \Nl^2+\right. \nonumber\\
		& \left.  2 \s_{0,1,1,0} \s_{1,0,0,1} \Nl^2+2 \s_{0,1,0,1} \s_{1,0,1,0} \Nl^2-4 \s_{0,1,0,0} \s_{1,0,1,1} \Nl^2-\right.\nonumber\\ 
		& \left. 4 \s_{0,0,1,1} \s_{1,1,0,0} \Nl^2+ 6 \s_{0,1,0,0} \s_{1,0,0,0} \s_{1,1,0,0} \Nl^2-4 \s_{0,0,1,0} \s_{1,1,0,1} \Nl^2-2 \s_{1,1,1,1} \Nl^2+\right.\nonumber\\
		& \left. 4 \s_{1,0,0,0} \s_{1,2,0,0} \Nl^2+\s_{0,1,0,0}^2 \s_{2,0,0,0} \Nl^2-\s_{0,2,0,0} \s_{2,0,0,0} \Nl^2+ 4 \s_{0,1,0,0} \s_{2,1,0,0} \Nl^2+ \right.\nonumber\\
		& \left. \s_{2,2,0,0} \Nl^2+2 \s_{0,0,1,1}^2 \Nl-4 \s_{0,1,0,0}^2 \s_{1,0,0,0}^2 \Nl-\s_{0,2,0,0} \s_{1,0,0,0}^2 \Nl+ \right.\nonumber\\
		& \left. 2 \s_{1,1,0,0}^2 \Nl-2 \s_{0,0,1,0} \s_{0,0,1,2} \Nl+(\Nl-1) \s_{0,0,0,2} \left(\s_{0,0,1,0}^2-\s_{0,0,2,0}\right) \Nl+ \right.\nonumber\\
		& \left. 8 \s_{0,0,1,1} \s_{0,1,0,0} \s_{1,0,0,0} \Nl+2 \s_{0,0,1,0} \s_{0,1,0,1} \s_{1,0,0,0} \Nl+ 2 \s_{0,1,1,1} \s_{1,0,0,0} \Nl+\right.\nonumber\\
		& \left. 2 \s_{0,0,1,0} \s_{0,1,0,0} \s_{1,0,0,1} \Nl-2 \s_{0,1,1,0} \s_{1,0,0,1} \Nl-2 \s_{0,1,0,1} \s_{1,0,1,0} \Nl+ \right.\nonumber\\
		& \left. 2 \s_{0,1,0,0} \s_{1,0,1,1} \Nl-2 \s_{0,0,1,1} \s_{1,1,0,0} \Nl-10 \s_{0,1,0,0} \s_{1,0,0,0} \s_{1,1,0,0} \Nl+ \right.\nonumber\\
		& \left. 2 \s_{0,0,1,0} \s_{1,1,0,1} \Nl-2 \s_{1,0,0,0} \s_{1,2,0,0} \Nl-\s_{0,1,0,0}^2 \s_{2,0,0,0} \Nl+ \right.\\
		& \left. \s_{0,2,0,0} \s_{2,0,0,0} \Nl-2 \s_{0,1,0,0} \s_{2,1,0,0} \Nl+6 \s_{0,1,0,0}^2 \s_{1,0,0,0}^2+ \right.\nonumber\\
		& \left. \s_{0,0,0,1}^2 \left((6-4 \Nl) \s_{0,0,1,0}^2+(\Nl-1) \Nl \s_{0,0,2,0}\right)-2 \s_{0,0,0,1} \right.\nonumber\\
		& \left. \left(\s_{0,0,1,0} \left(\Nl (5-3 \Nl) \s_{0,0,1,1}+2 (3-2 \Nl) \s_{0,1,0,0} \s_{1,0,0,0}+ 2 (\Nl-2) \Nl \s_{1,1,0,0}\right)+ \right. \right.\nonumber\\
		& \left. \left. (\Nl-1) \Nl \left((\Nl-1) \s_{0,0,2,1}+\s_{0,1,1,0} \s_{1,0,0,0}+\right. \right.\right.\nonumber\\
		&  \left. \left. \left. \s_{0,1,0,0} \s_{1,0,1,0}-\Nl \s_{1,1,1,0}+\s_{1,1,1,0}\right)\right)\right) \nonumber.
\end{align}
\textcolor{black}{Here, the sample size on each mesh level $l$ must be $\Nl>3$.}

\subsubsection{Computational Complexity}

If $\mathcal{C}(\Zl)\in\mathbb{R}_{\geq 0}$ represents the computational cost of evaluating one MC sample of $u^{h_l}(x, \event)$ and $u^{h_{l-1}}(x, \event)$ in the difference term $\Zl$, then the overall cost of the MLMC covariance estimator $\honeoneML$ is given by:
\begin{equation}\label{Eq:costML}
	\mathcal{C}(\honeoneML) := \sum_{\meshlevel=0}^{\LL} \Nl \mathcal{C}(\Zl).
\end{equation}

Subsequently, the optimum number of samples $\Nl$ on each mesh level $l$ is determined to minimize the total computational cost $\mathcal{C}(\honeoneML)$ under the constraint $\sum_{\meshlevel=0}^{\numofmesh} \varMC(\Zl) < \epsilon^2/2$. 
By re-expressing $\varMC(\Zl) = \vloneone/\Nl$, which essentially is $\vloneone = \Nl \, \varMC(\Zl)$, this translates into a constrained optimization problem, where the cost function
\begin{equation}\label{Eq:costfn}
	f(\Nl) = \text{arg} \min_{\optSam} \sum_{\meshlevel=0}^{\LL} \left(\optSam \mathcal{C}(\Zl) + \Lagrange \frac{\vloneone}{\Nl}\right),
\end{equation}
is minimized. Therefore, we can obtain the optimal number of samples $\Nl$ as:
\begin{equation}\label{eq:Nl}
	\optSam = \Lagrange\left({\frac{\vloneone}{\mathcal{C}(\Zl)}}\right)^{\frac{1}{2}},
\end{equation}
where $\Lagrange$ is the Lagrange multiplier defined by
\begin{equation}\label{Eq:tau}
	\Lagrange = {\frac{2}{\epssq}} \sum_{\meshlevel=0}^{\numofmesh} \left({{\vloneone\mathcal{C}(\Zl)}}\right) ^{\frac{1}{2}}.
\end{equation}

In addition, one can strategically choose the mesh levels $\left\lbrace \meshlevel=0,1,2,...,\LL \right\rbrace$ in an optimal manner, either following a geometric or non-geometric sequence, to minimize a specified deterministic error. For example, in PDE-based applications, the selection of the initial coarse mesh $\meshlevel=0$ is contingent on the regularity of the solution $u(x)$ in \feq{abs_main_eq} \cite{Cliffe2011}. Subsequently, the choice of finer mesh grids is guided by an a priori mesh convergence analysis, where the finest mesh level $\meshlevel = L$ can either be set in advance or adaptively determined \cite{giles_2015} during the computation of the MLMC estimation. For a deeper exploration of the optimal strategies for mesh hierarchy selection, refer to \cite{haji-ali_optimization_2016}.

The overall computational complexity of the h-statistics-based MLMC covariance estimator mirrors that of the classical MLMC procedure \rvsi{of mean estimation in \cite{Cliffe2011,giles_2015}. Furthermore, the notion of computational complexity of an arbitrary multilevel statistical estimator is detailed in \cite[Theorem 2.1]{mycek_multilevel_2019}, which also derives a special case for the multilevel covariance estimation in Theorem 2.5.} {Therefore, we consider the following proposition:}
\begin{proposition}\label{prop:mlcov}
	Consider the positive constants \(c_\alpha\), \(c_\beta\), \(c_\gamma\), \(\convorder\), \(\convbeta\), and \(\convgamma\) with the constraint \(\convorder \geq \frac{1}{2}\min(\convbeta, \convgamma)\). We assume the following error bounds:
	\begin{equation}\label{eq:errorcov}
		\begin{aligned}
			(i)&\ \text{the deterministic error decays as}\ |\mu_{1,1}(u^{h_l}_i, u^{h_l}_j)- \mu_{1,1}\uij| \leq c_\alpha {h_l}^{\alpha},\\
			(ii)&\ \text{the variance decays as}\ \vloneone \leq c_\beta h_l^{\convbeta}, \\
			(iii)&\ \text{the computational cost to determine one sample of} \ \Zl \ \text{is given by} \\ 
			&\quad \mathcal{C}(\Zl)  \leq c_\convgamma {h_l}^{-\convgamma}.
		\end{aligned}
	\end{equation}
	Under these conditions, there exists a positive constant \(c\) such that, for any \(0 < \epsilon < e^{-1}\), the multilevel Monte Carlo (MLMC) mean square error estimator satisfies \(\mathrm{MSE}(\honeoneML) < \epsilon^2\), and the total computational cost is bounded as follows:
	\begin{equation}\label{Eq:costcomplex}
		\mathcal{C}(\honeoneML) \leq
		\begin{cases}
			c\,{\epsilon}^{-2}, & \convbeta > \convgamma,\\
			c\,{\epsilon}^{-2}\rvso{(\log{\epsilon})^2}, & \convbeta = \convgamma,\\
			c\,{\epsilon}^{-2-(\convgamma-\convbeta)/\convorder}, & \convbeta < \convgamma.
		\end{cases}
	\end{equation}
\end{proposition}

By considering the values of $\convbeta$ and $\convgamma$, we can determine the primary cost driver within the mesh grid sequence. If $\convbeta$ is greater than $\convgamma$, then the maximum cost is influenced by the coarsest level. Conversely, if $\convbeta$ is less than $\convgamma$, the finest level dictates the predominant cost. When $\convbeta$ equals $\convgamma$, the cost is generally evenly distributed across all levels

Furthermore, it is evident from the first relation in \feq{eq:errorcov} that, as \(l\rightarrow \infty\),
\[
|\mu_{1,1}(u^{h_l}_i, u^{h_l}_j) - \mu_{1,1}\uij| \rightarrow 0.
\]
However, \(\mu_{1,1}\uij\) is analytical and unknown. Therefore, the deterministic error is defined using the triangle inequality as
\begin{equation}\label{Eq:deterministic_error}
	|\mu_{1,1}(u^{h_l}_i, u^{h_l}_j) - \mu_{1,1}(u^{h_{l-1}}_i, u^{h_{l-1}}_j)| \leq c_\alpha {h_l}^{\alpha}.
\end{equation}

\section{Application to Poisson's equation: 1D steady state heat equation} \label{Sec:application}

\noindent \textbf{Deterministic Setting:} Consider a physical body defined in the domain $\spatialdomain\subset\mathbb{R}$, within a one-dimensional Euclidean space $\mathbb{R}$, with a smooth Lipschitz boundary $\Gamma$. The objective is to determine the temperature $\uu({x}) \in \mathbb{R}$ (belonging to the Hilbert space $\solspace$) at a spatial point ${x}\in \spatialdomain$. This solution should satisfy the governing equations describing the steady-state heat conductivity:
\begin{equation}\label{Eq:deterministic-cond}
	\begin{aligned}
		-\nabla\cdot(\kappa(x) \nabla \uu(x)) &= \heatsource(x)\quad \forall x\in\spatialdomain, \\
		\uu(x) &= {\uu}_0,  \quad \forall x \in \Gamma_D, \\
		\kappa(x) \nabla \uu(x)\cdot{n}(x) &= {g}(x), \quad \forall x \in \Gamma_N,
	\end{aligned}
\end{equation}
where $\heatsource(x) \in \mathbb{R}$ represents the heat sink and source. In addition, ${g}(x) \in \mathbb{R}$ and $\uu_0 \in \mathbb{R}$ denote the surface heat flux on the Neumann boundary $\Gamma_N \subset \Gamma$ and the homogeneous temperature on the Dirichlet boundary $\Gamma_D \subseteq \Gamma$, respectively. The vector $n(x) \in \Euclideanspace$ signifies the outward unit normal to the boundary $\Gamma_N$.
Furthermore, we assume that $\Gamma = \Gamma_D \cup \Gamma_N$ and $\Gamma_D \cap \Gamma_N = 0$.

In this context, $\kappa(x)\in\mathbb{R}^+$ represents the thermal conductivity field. As Poisson's equation encompasses a broader class of elliptic boundary value problems, describing physical phenomena such as electrostatic potential and chemical diffusion, the material coefficient $\kappa(x)$ may also be considered to characterize intrinsic physics-based properties, such as electric permittivity and chemical diffusivity.

\noindent \textbf{Uncertainty Modelling:} In accordance with \fsec{UncertaintyQuantification}, we specifically model the coefficient $\kappa$ as uncertain. For simplicity, this study only considers a spatially constant thermal conductivity such that $\kappa(x) \equiv \kappa\in\mathbb{R}^+$.
Thus, $\kappa(x)$ is depicted as a positive-definite random variable
\begin{equation}\label{Eq:Tmapping}
	\kappa(\event): \samplespace\rightarrow\mathbb{R}^+
\end{equation}
over probability space $(\samplespace,\sigmaalgebra,\probability)$.

For instance, in this study, we model $\kappa(\event)$ as a log-normal random variable
\begin{equation}\label{Eq:lognormal}
	\kappa(\event) = \exp(\psi(\event))
\end{equation}
on $\mathbb{R}^+$. Here, $\psi(\event)$ represents a Gaussian random variable,
\begin{equation}
	\psi(\event) = \log(\kappa(\event)) \sim \normaldist(\mu,\sigma)
\end{equation}
over $\mathbb{R}$, which is characterized by the mean $\mu\in\mathbb{R}$ and standard deviation $\sigma\in\PosRealStrict$.

The desired parameters, namely the mean $\mu_{\kappa}$ and standard deviation $\sigma_{\kappa}$ of the log-normal variable $\kappa(\event)$ in \feq{Eq:lognormal}, are related to the corresponding Gaussian random variable $\psi(\event)$ as established in \cite{pinsky_1_2011}:
\begin{equation}
	\mu = \ln\left(\frac{\mu_{\kappa}^2}{\sqrt{\mu_{\kappa}^2 + \sigma_{\kappa}^2}}\right)  \quad \mathrm{and} \quad \sigma =\sqrt{\ln \left(1 + \frac{\sigma_{\kappa}^2}{\mu_{\kappa}^2}\right)}.
\end{equation}

\noindent \textbf{Stochastic Setting:}
By reformulating the deterministic model in \feq{Eq:deterministic-cond} as a stochastic model,
\begin{equation}\label{Eq:stochastic-thm}
	\begin{aligned}
		-\nabla\cdot(\kappa(\event) \nabla \uu(x,\event)) &= \heatsource(x)\quad \forall x\in\spatialdomain, \event\in\samplespace, \\
		\uu(x,\event) &= {\uu}_0,  \quad \forall x \in \Gamma_D, \\
		\kappa(\event) \nabla \uu(x,\event) \cdot n(x) &= {g}(x), \quad \forall x \in \Gamma_N,
	\end{aligned}
\end{equation}
which must be understood in a weak sense, both spatially and stochastically.
The objective is to determine the random temperature field $\temperature({x},\event):\spatialdomain\times\samplespace\rightarrow\mathbb{R}$ by assuming deterministic boundary conditions and a heat source.

As explained in \fsec{UncertaintyQuantification}, given the weak form of \feq{Eq:stochastic-thm}, we employ the finite-element method to seek an approximate solution $u^h({x},\event)$ within a finite-dimensional subspace $\mathcal{U}^h\subset\solspace$ of the solution space $\solspace$, where $h$ denotes the discretization parameter, which is an indicator of element size.
Subsequently, this study aims to determine the sample covariance of the solution $u^h({x},\event)$ using Monte Carlo and multilevel Monte Carlo estimators, as defined in \fsecs{Sec:MC}{Sec:MLMC}.

\section{Numerical results}\label{NumericalResults}

\subsection{Specifications}

\begin{figure}[ht]
	\centering
	\def\svgwidth{14.3cm}
\begingroup%
  \makeatletter%
  \providecommand\color[2][]{%
    \errmessage{(Inkscape) Color is used for the text in Inkscape, but the package 'color.sty' is not loaded}%
    \renewcommand\color[2][]{}%
  }%
  \providecommand\transparent[1]{%
    \errmessage{(Inkscape) Transparency is used (non-zero) for the text in Inkscape, but the package 'transparent.sty' is not loaded}%
    \renewcommand\transparent[1]{}%
  }%
  \providecommand\rotatebox[2]{#2}%
  \newcommand*\fsize{\dimexpr\f@size pt\relax}%
  \newcommand*\lineheight[1]{\fontsize{\fsize}{#1\fsize}\selectfont}%
  \ifx\svgwidth\undefined%
    \setlength{\unitlength}{465.48449419bp}%
    \ifx\svgscale\undefined%
      \relax%
    \else%
      \setlength{\unitlength}{\unitlength * \real{\svgscale}}%
    \fi%
  \else%
    \setlength{\unitlength}{\svgwidth}%
  \fi%
  \global\let\svgwidth\undefined%
  \global\let\svgscale\undefined%
  \makeatother%
	\footnotesize
  \begin{picture}(1,0.09533146)%
    \lineheight{1}%
    \setlength\tabcolsep{0pt}%
    \put(0,0){\includegraphics[width=\unitlength,page=1]{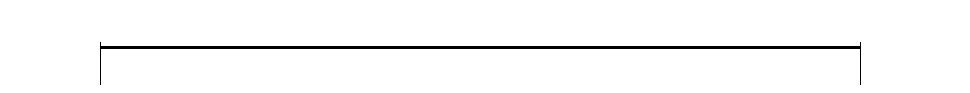}}%
    \put(0.08277822,0.05967148){\color[rgb]{0,0,0}\makebox(0,0)[lt]{\lineheight{1.25}\smash{\begin{tabular}[t]{l}$u(0) = 273 \, \mathrm{K}$\end{tabular}}}}%
    \put(0.8015526,0.05833718){\color[rgb]{0,0,0}\makebox(0,0)[lt]{\lineheight{1.25}\smash{\begin{tabular}[t]{l}$u(1) = 273\, \mathrm{K}$\end{tabular}}}}%
    \put(0.45003441,0.00029491){\color[rgb]{0,0,0}\makebox(0,0)[lt]{\lineheight{1.25}\smash{\begin{tabular}[t]{l}1 m\end{tabular}}}}%
    \put(0,0){\includegraphics[width=\unitlength,page=2]{1D_heat.pdf}}%
    \put(0.41053118,0.05833718){\color[rgb]{0,0,0}\makebox(0,0)[lt]{\lineheight{1.25}\smash{\begin{tabular}[t]{l}$g(x) = 5$ W/m\end{tabular}}}}%
  \end{picture}%
\endgroup%

	\caption{Geometry and boundary conditions}
	\label{Fig:BC1D}
\end{figure}

We consider a one-dimensional geometry with a body width of approximately \(1 \, \text{m}\), as illustrated in \ffig{Fig:BC1D}. A uniformly distributed surface heat flux of \(5 \, \text{W/m}\) is applied over the entire length, and the geometry's ends are maintained at a fixed temperature of \(273 \, \text{K}\).

\begin{figure}[ht]
	\centering
\begin{tikzpicture}

\definecolor{crimson2143940}{RGB}{214,39,40}
\definecolor{darkgray176}{RGB}{176,176,176}
\definecolor{darkorange25512714}{RGB}{255,127,14}
\definecolor{forestgreen4416044}{RGB}{44,160,44}
\definecolor{steelblue31119180}{RGB}{31,119,180}
\begin{axis}[
width=14cm,
height=12cm,		
hide x axis,
hide y axis,
tick align=outside,
tick pos=left,
x grid style={darkgray176},
xmin=-0.05, xmax=1.05,
xtick style={color=black},
y grid style={darkgray176},
ymin=0.5, ymax=1.2,
ytick style={color=black}
]
\addplot [line width=2pt, mark=*, mark size=1.0, mark options={solid}]
table {%
0 1
0.1 1
0.2 1
0.3 1
0.4 1
0.5 1
0.6 1
0.7 1
0.8 1
0.9 1
1 1
};
\addplot [line width=2pt, mark=*, mark size=1.0, mark options={solid}]
table {%
0 0.9
0.05 0.9
0.1 0.9
0.15 0.9
0.2 0.9
0.25 0.9
0.3 0.9
0.35 0.9
0.4 0.9
0.45 0.9
0.5 0.9
0.55 0.9
0.6 0.9
0.65 0.9
0.7 0.9
0.75 0.9
0.8 0.9
0.85 0.9
0.9 0.9
0.95 0.9
1 0.9
};
\addplot [line width=2pt, mark=*, mark size=1.0, mark options={solid}]
table {%
0 0.8
0.025 0.8
0.05 0.8
0.075 0.8
0.1 0.8
0.125 0.8
0.15 0.8
0.175 0.8
0.2 0.8
0.225 0.8
0.25 0.8
0.275 0.8
0.3 0.8
0.325 0.8
0.35 0.8
0.375 0.8
0.4 0.8
0.425 0.8
0.45 0.8
0.475 0.8
0.5 0.8
0.525 0.8
0.55 0.8
0.575 0.8
0.6 0.8
0.625 0.8
0.65 0.8
0.675 0.8
0.7 0.8
0.725 0.8
0.75 0.8
0.775 0.8
0.8 0.8
0.825 0.8
0.85 0.8
0.875 0.8
0.9 0.8
0.925 0.8
0.95 0.8
0.975 0.8
1 0.8
};
\addplot [line width=2pt, mark=*, mark size=1.0, mark options={solid}]
table {%
0 0.7
0.0125 0.7
0.025 0.7
0.0375 0.7
0.05 0.7
0.0625 0.7
0.075 0.7
0.0875 0.7
0.1 0.7
0.1125 0.7
0.125 0.7
0.1375 0.7
0.15 0.7
0.1625 0.7
0.175 0.7
0.1875 0.7
0.2 0.7
0.2125 0.7
0.225 0.7
0.2375 0.7
0.25 0.7
0.2625 0.7
0.275 0.7
0.2875 0.7
0.3 0.7
0.3125 0.7
0.325 0.7
0.3375 0.7
0.35 0.7
0.3625 0.7
0.375 0.7
0.3875 0.7
0.4 0.7
0.4125 0.7
0.425 0.7
0.4375 0.7
0.45 0.7
0.4625 0.7
0.475 0.7
0.4875 0.7
0.5 0.7
0.5125 0.7
0.525 0.7
0.5375 0.7
0.55 0.7
0.5625 0.7
0.575 0.7
0.5875 0.7
0.6 0.7
0.6125 0.7
0.625 0.7
0.6375 0.7
0.65 0.7
0.6625 0.7
0.675 0.7
0.6875 0.7
0.7 0.7
0.7125 0.7
0.725 0.7
0.7375 0.7
0.75 0.7
0.7625 0.7
0.775 0.7
0.7875 0.7
0.8 0.7
0.8125 0.7
0.825 0.7
0.8375 0.7
0.85 0.7
0.8625 0.7
0.875 0.7
0.8875 0.7
0.9 0.7
0.9125 0.7
0.925 0.7
0.9375 0.7
0.95 0.7
0.9625 0.7
0.975 0.7
0.9875 0.7
1 0.7
};
\draw (axis cs:0,1.03) node[
  scale=0.7,
  anchor=base west,
  text=black,
  rotate=0.0
]{Mesh level=0, $E_0$=10, $M_0$=11};
\draw (axis cs:0,0.93) node[
  scale=0.7,
  anchor=base west,
  text=black,
  rotate=0.0
]{Mesh level=1, $E_1$=20, $M_1$=21};
\draw (axis cs:0,0.83) node[
  scale=0.7,
  anchor=base west,
  text=black,
  rotate=0.0
]{Mesh level=2, $E_2$=40, $M_2$=41};
\draw (axis cs:0,0.73) node[
  scale=0.7,
  anchor=base west,
  text=black,
  rotate=0.0
]{Mesh level=3, $E_3$=80, $M_3$=81};
\end{axis}

\end{tikzpicture}
	\caption{Four levels of nested mesh grids}
	\label{Fig:mesh1D}
\end{figure}
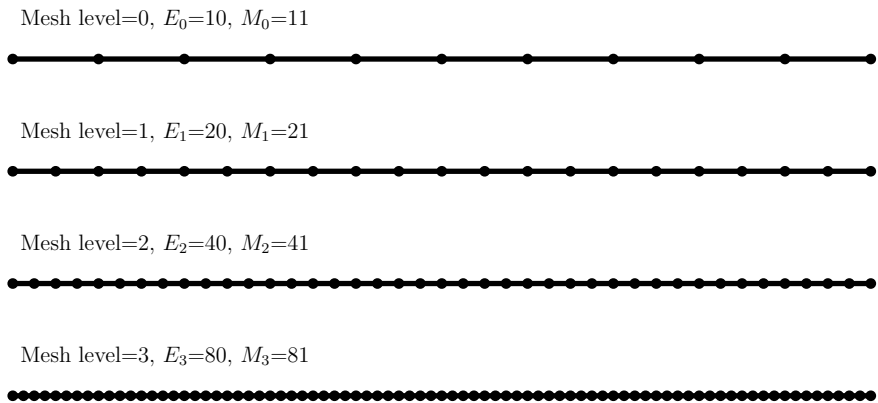

The spatial discretization in this study is based on the Finite Element Method and utilizes two-node 1D linear elements. To implement MLMC, we adopt a sequence of four nested meshes, each with an element size of \(h_{\meshlevel-1} = 2h_\meshlevel\). For specific mesh specifications, please refer to \ffig{Fig:mesh1D}. In this context, \(E_l\) represents the number of elements, and \(M_l\) indicates the total spatial nodes at each mesh level \(l\). Identical boundary conditions are consistently applied across all levels of the mesh.

As described in \fsec{Sec:application}, we model the thermal conductivity \(\kappa(\event)\) as a log-normal random variable with parameters: \(\mu_{\kappa} = 0.1 \, \text{W/mK}\) and standard deviation \(\sigma_{\kappa} = 0.03 \, \text{W/mK}\). The objective of this study is to determine the sample covariance of the stochastic solution \(u^h({x},\event)\) at \(L=3\) using h-statistics-based Monte Carlo (\(\honeoneMC\), see \feq{Eq:h11}), classical Multilevel Monte Carlo (\(\covML\), defined in \feq{Eq:covMLMC}), and h-statistics-based Multilevel Monte Carlo estimators (\(\honeoneML\), refer to \feq{Eq:covMLMChstat}).

The analytical response \(u(x,\event)\) in \feq{Eq:stochastic-thm} is spatially discretized at \(M_l\) spatial nodes at each mesh level \(l\); clearly, \(M_0 < M_1 < M_2 < M_3\). Consequently, to determine the sample covariance matrix (in the form shown in \feq{Eq:cov_analytical_vec}) on \(L=3\), i.e., \(\honeoneML(\vek{U}^{h_L}) \in \PosSymMatSpace(\nodes_L)\), the corresponding MC estimators in \feq{Eq:covMLMChstat} on mesh levels \(l < L\) are interpolated with respect to the finest level \(L\).

Furthermore, to facilitate handling fully filled sample covariance matrices, we utilize half-vectorization for a positive-definite matrix \(\vek{H} \in \PosSymMatSpace(\nodes)\):
\begin{equation}
	\vech(\vek{H}) := [H_{1,1}, \ldots, H_{M, 1}, H_{2,2}, \ldots, H_{M,2}, \ldots, H_{M-1,M-1}, H_{M,M-1}, H_{M,M}]^{\mathrm{T}}.
\end{equation}
Here, the vector \(\vech(\vek{H}) \in \mathbb{R}^{M(M+1)/2}\).

Assuming the optimal finest mesh grid \(L=3\) is known, this implies that the square of the deterministic error is less than \(\epsilon^2/2\) as indicated in \feqss{Eq:MSEhstat}{Eq:MSEcovML}{Eq:MSEhstatML}. Consequently, the objective of this paper is to analyse and ensure that only the stochastic accuracy is less than \(\epsilon^2/2\).

\subsection{Screening}

An a priori performance analysis of h-statistics based MLMC covariance estimator is conducted using a screening test. In this test, a fixed number of 50 samples is examined over four levels of meshes.

By considering the ${L}_\infty$ or $\max$ norm over the estimated and half-vectorized error bounds in \feq{Eq:deterministic_error}, we obtain:
\begin{equation}\label{}
		\max|\vech({\honeoneMC}(\vek{U}^{h_l}) - {\honeoneMC(\vek{U}^{h_{l-1}}}))| \leq c_\alpha {h_l}^{\alpha},
\end{equation}
which is equivalent to
\begin{equation}\label{Eq:deterr}
	\max|\vech(\Zl)| \leq c_\alpha {h_l}^{\alpha}.
\end{equation}
Refer to \feq{Eq:Zlhstat} for the definition of $\Zl$. Similarly, the second relation in \feq{eq:errorcov} is expressed as:
\begin{equation}\label{Eq:stoerr}
	 \max(\vech(\vloneone)) \leq c_\beta h_l^{\convbeta}.
\end{equation}

\rvsn{To analyze the performance difference between the classical and h-statistics MLMC covariance estimator, we conduct an additional screening test of the conventional estimator. Subsequently, referring to \feq{eq:appsamperrul}, we define
$$\varMC(\Yl) = \widehat{\mathbb{V}}_{l,\mathrm{cov}}^{\mathrm{MC}} / \Nl;$$ in other words, $\widehat{\mathbb{V}}_{l,\mathrm{cov}}^{\mathrm{MC}} = \Nl \varMC(\Yl)$. Thus, similar to \feq{Eq:stoerr}, we obtain
\begin{equation}\label{Eq:stoerrul}
	\max(\vech(\widehat{\mathbb{V}}_{l,\mathrm{cov}}^{\mathrm{MC}})) \leq c_\beta^* h_l^{\convbeta^*},
\end{equation}
where $c_\beta^*$ and $\beta^*$ are positive constants.
It is important to note that the deterministic decay and the propagation of computational cost, as presented in \feq{Eq:deterr} and the third expression in \feq{eq:errorcov}, respectively, remain identical for the classical and the h-statistics estimator.}

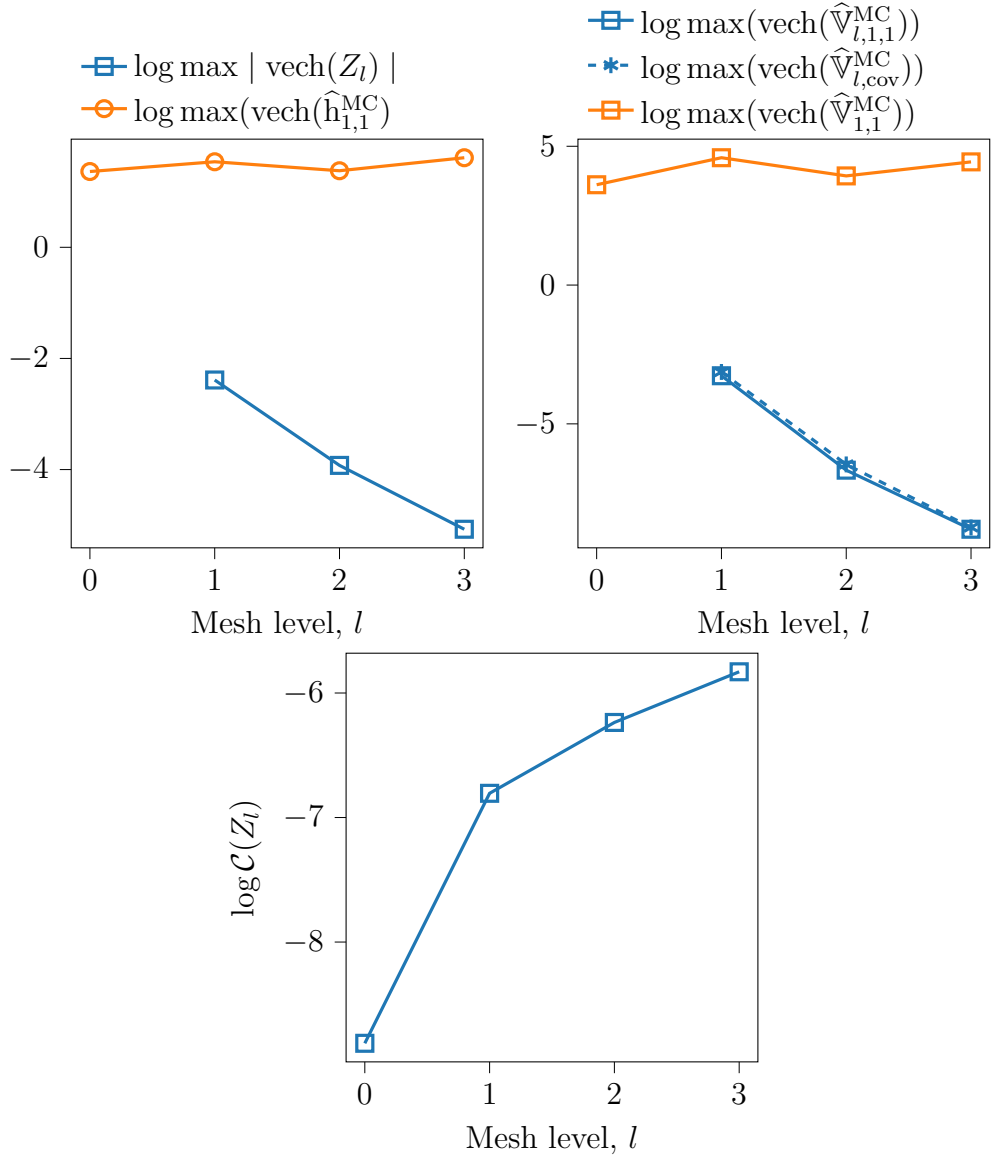
\begin{figure}[ht!]
	\centering
\begin{tikzpicture}
\definecolor{darkgray176}{RGB}{176,176,176}
\definecolor{darkorange25512714}{RGB}{255,127,14}
\definecolor{lightgray204}{RGB}{204,204,204}
\definecolor{steelblue31119180}{RGB}{31,119,180}

\begin{axis}[
width=7cm,
height=7cm,	
legend cell align={left},
legend style={
  fill opacity=0.8,
  draw opacity=1,
  text opacity=1,
  at={(0.03,0.03)},
  anchor=south west,
  draw=lightgray204
},
tick align=outside,
tick pos=left,
x grid style={darkgray176},
xlabel={Mesh level, $l$},
xmin=-0.15, xmax=3.15,
xtick style={color=black},
y grid style={darkgray176},
ymin=-5.40813373198015, ymax=1.94597013326333,
ytick style={color=black},
axis background/.style={fill=white},
legend style={at={(0,1.24)}, anchor=north west, legend cell align=left, align=left, draw=none}
]
\addplot [semithick, steelblue31119180, line width=1.2pt, mark=square, mark size=3, mark options={solid}]
table {%
1 -2.38801749903285
2 -3.92466052353714
3 -5.07385628355999
};
\addlegendentry{$\log \mathrm{max}\mid \vech(\Zl)\mid$}
\addplot [semithick, darkorange25512714, line width=1.2pt, mark=o, mark size=3, mark options={solid}]
table {%
0 1.36347953089217
1 1.53911911714142
2 1.37741382162109
3 1.61169268484317
};
\addlegendentry{$\log \mathrm{max}(\vech(\honeoneMC)$}
\end{axis}

\end{tikzpicture}
	\centering
\begin{tikzpicture}
\definecolor{darkgray176}{RGB}{176,176,176}
\definecolor{darkorange25512714}{RGB}{255,127,14}
\definecolor{lightgray204}{RGB}{204,204,204}
\definecolor{steelblue31119180}{RGB}{31,119,180}

\begin{axis}[
width=7cm,
height=7cm,	
legend cell align={left},
legend style={
  fill opacity=0.8,
  draw opacity=1,
  text opacity=1,
  at={(0.03,0.03)},
  anchor=south west,
  draw=lightgray204
},
tick align=outside,
tick pos=left,
x grid style={darkgray176},
xlabel={Mesh level, $l$},
xmin=-0.15, xmax=3.15,
xtick style={color=black},
y grid style={darkgray176},
ymin=-9.46586851872904, ymax=5.25669789312267,
ytick style={color=black},
axis background/.style={fill=white},
legend style={at={(0,1.36)}, anchor=north west, legend cell align=left, align=left, draw=none}
]
\addplot [semithick, steelblue31119180, line width=1.2pt, mark=square, mark size=3, mark options={solid}]
table {%
1 -3.26683512036009
2 -6.66842924870469
3 -8.79666095455396
};
\addlegendentry{$\log \mathrm{max}(\vech(\vloneone))$}
\addplot [semithick, steelblue31119180, line width=1.2pt, , dashed, mark=asterisk, mark size=3, mark options={solid}]
table {%
	1 -3.130389652333893
	2 -6.461027698552701
	3 -8.731723932613678
};
\addlegendentry{$\log \mathrm{max}(\vech( \widehat{\mathbb{V}}_{l,\mathrm{cov}}^{\mathrm{MC}} ))$}
\addplot [semithick, darkorange25512714, line width=1.2pt, mark=square, mark size=3, mark options={solid}]
table {%
0 3.61415189
1 4.58749032894759
2 3.93220113696243
3 4.43772735685716
};
\addlegendentry{$\log \mathrm{max}(\vech(\voneoneMC))$}
\end{axis}

\end{tikzpicture}
	\centering
\begin{tikzpicture}
\definecolor{darkgray176}{RGB}{176,176,176}
\definecolor{steelblue31119180}{RGB}{31,119,180}

\begin{axis}[
width=7cm,
height=7cm,	
tick align=outside,
tick pos=left,
x grid style={darkgray176},
xlabel={Mesh level, $l$},
xmin=-0.15, xmax=3.15,
xtick style={color=black},
y grid style={darkgray176},
ylabel={ $ \log  \mathcal{C}(\Zl) $},
ymin=-8.96225685874195, ymax=-5.68038038424494,
ytick style={color=black}
]
\addplot [semithick, steelblue31119180, line width=1.2pt, mark=square, mark size=3, mark options={solid}]
table {%
0 -8.81308065535572
1 -6.80552966626566
2 -6.23649379775459
3 -5.82955658763116
};
\end{axis}

\end{tikzpicture}
	\caption{Screening results of MLMC covariance estimators $\honeoneML$ \rvsn{and $\covML$}}
	\label{Fig:screening}
\end{figure}

An overview of the corresponding results is presented in \ffig{Fig:screening}. The top left and right plots illustrate the behaviour of the logarithm of quantities defined in \feq{Eq:deterr} \rvsn{and \feqs{Eq:stoerr}{Eq:stoerrul},} respectively. The former diagram demonstrates deterministic decay of $\max|\vech(\Zl)|$, while the latter showcases stochastic convergences \rvsn{of h-statistics and classical estimator i.e., $ \max(\vech(\vloneone))$ and $\mathrm{max}(\vech(\widehat{\mathbb{V}}_{l,\mathrm{cov}}^{\mathrm{MC}} ))$, respectively. 
An interesting observation in this regard is that for all values of $\meshlevel>0$, $ \max(\vech(\vloneone))<\mathrm{max}(\vech(\widehat{\mathbb{V}}_{l,\mathrm{cov}}^{\mathrm{MC}} ))$. This suggests that, for a given accuracy \(\epsilon^2/2\), the number of samples $\Nl$ on each level $\meshlevel>0$, h-statistics-based MLMC covariance estimator will incur fewer samples than the classical one.}
Additionally, both plots depict the quantities $\mathrm{max}(\vech(\honeoneMC)$ and $\mathrm{max}(\vech(\voneoneMC))$, which remain approximately constant at all values of $\meshlevel$.


Moreover, the logarithmic computational time of running one sample $\mathcal{C}(\Zl)$ on each mesh level $\meshlevel$ is displayed at the bottom. These values are obtained by recording the timings for the 50 considered screening samples on a M1 Pro chip with 16GB of RAM and averaging them. As the results indicate, computation becomes more expensive as mesh refinement increases.

Finally, we determine the decay rates and constants corresponding to  \rvsn{\feqss{Eq:deterr}{Eq:stoerr}{Eq:stoerrul}}, and the third equation in \feq{eq:errorcov} by analyzing the slopes and y-intercepts of the respective logarithmic quantities in \ffig{Fig:screening}. The results are summarized in  \ftbl{Table:convergence_mean}.
It is evident that all constants are positive, and the condition $\convorder \geq \frac{1}{2}\text{min}(\convbeta,\convgamma)$ is satisfied, confirming the assumption made in \prop{prop:mlcov}. 
\rvsn{Furthermore, it is noteworthy that the decay rates associated with stochastic convergences by both covariance estimators are nearly equal, with $\beta \approx \beta^*$, while their respective constants satisfy $c_{\beta} < c_{\beta}^*$. This implies that for a given accuracy  \(\epsilon^2/2\), the computational cost of the h-statistics-based MLMC covariance estimator at each level $\meshlevel$ will be lower than that of the conventional MLMC estimator.}
Furthermore, with $\convbeta>\convgamma$, the computational complexity of MLMC estimates follows the first scenario outlined in \feq{Eq:costcomplex}.
\begin{table}[ht!]
	\centering
	\begin{tabular}{llllllll}
		\specialrule{1pt}{1pt}{1pt}
		 $\convorder$ & $\convbeta$ & \rvsn{$\convbeta^* $}  & $\convgamma$ & $c_{\alpha}$  & $c_{\beta}$ & \rvsn{$c_{\beta}^*$}  & $c_{\gamma}$  \\ \midrule
		 {$1.94$} &  {$3.99$} & \rvsn{4.04}  & {$1.37$} & 0.33 &  {$0.49$}  & \rvsn{0.60} &  {0.0002} \\
		\specialrule{1pt}{1pt}{1pt}
	\end{tabular}
	\caption{Convergence results of MLMC covariance estimators $\honeoneML$ \rvsn{and $\covML$}}
	\label{Table:convergence_mean}
\end{table} 

\subsection{Performance}

\subsubsection{Accuracy}
\begin{figure}[ht!]
		\centering
		\includegraphics[width=0.45\linewidth]{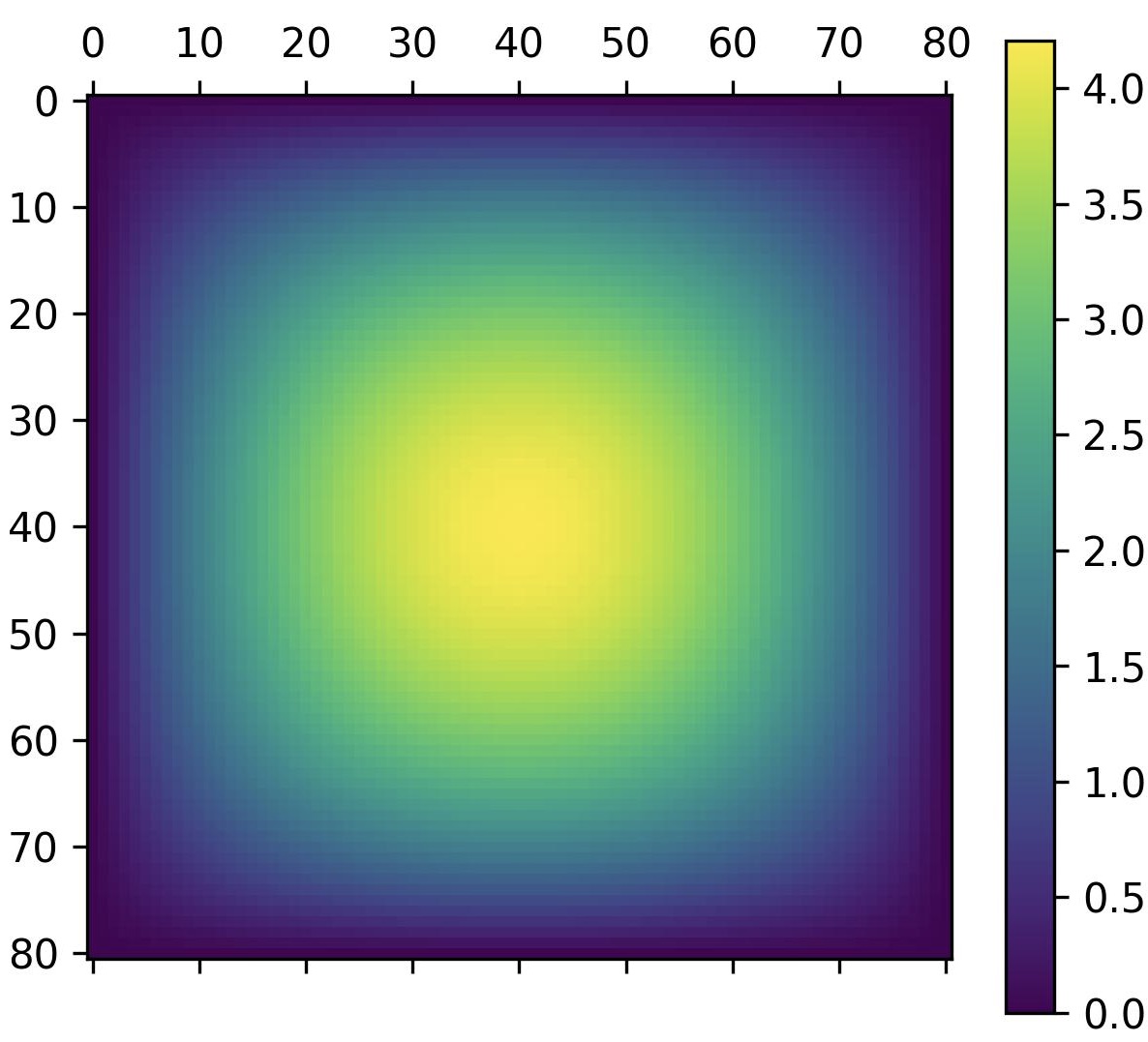} \hfill
		\includegraphics[width=0.45\linewidth]{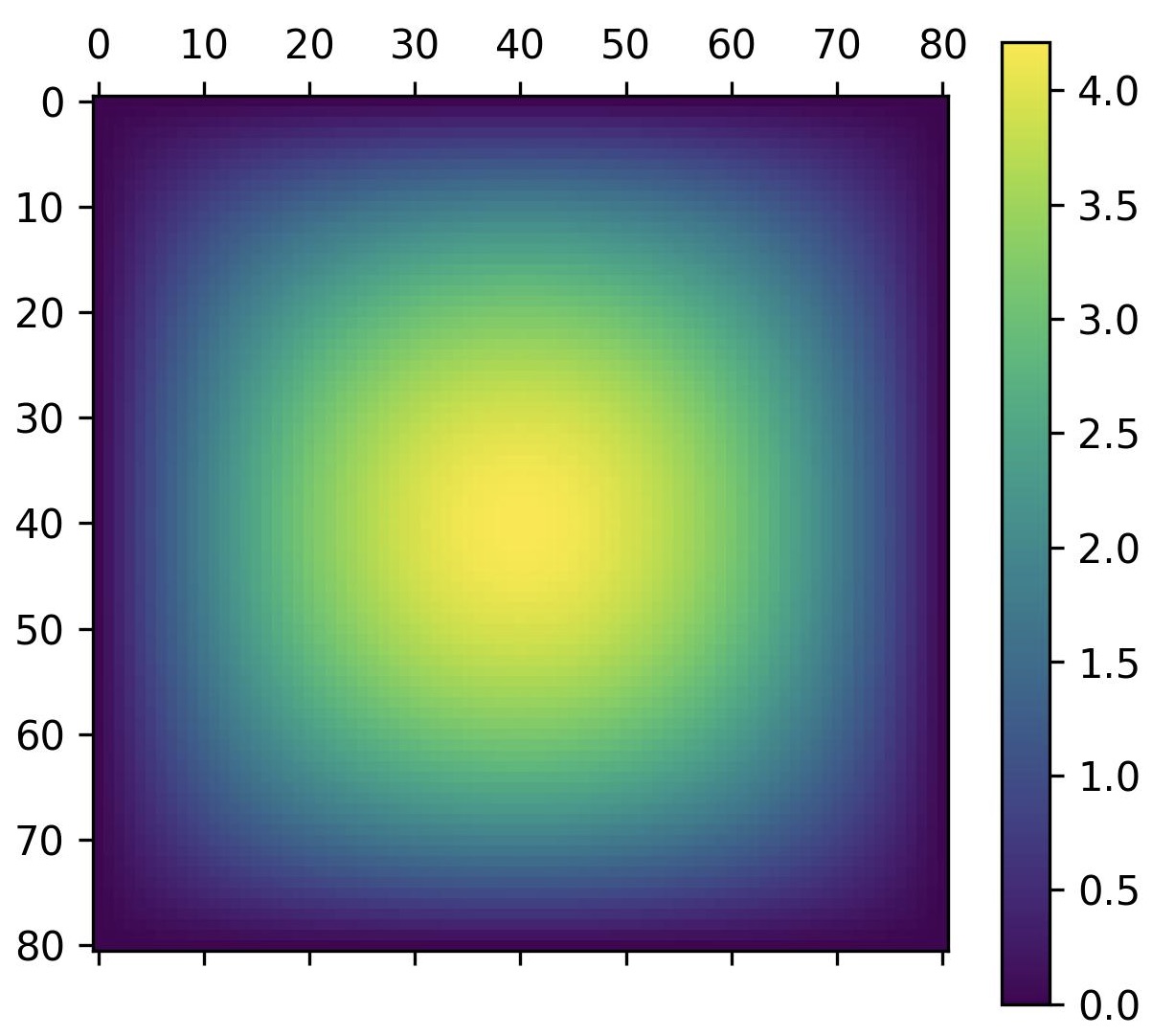} \\
		(a) $\honeoneML$ \hspace*{200pt}	(b) $\covML$ \\
		\includegraphics[width=0.45\linewidth]{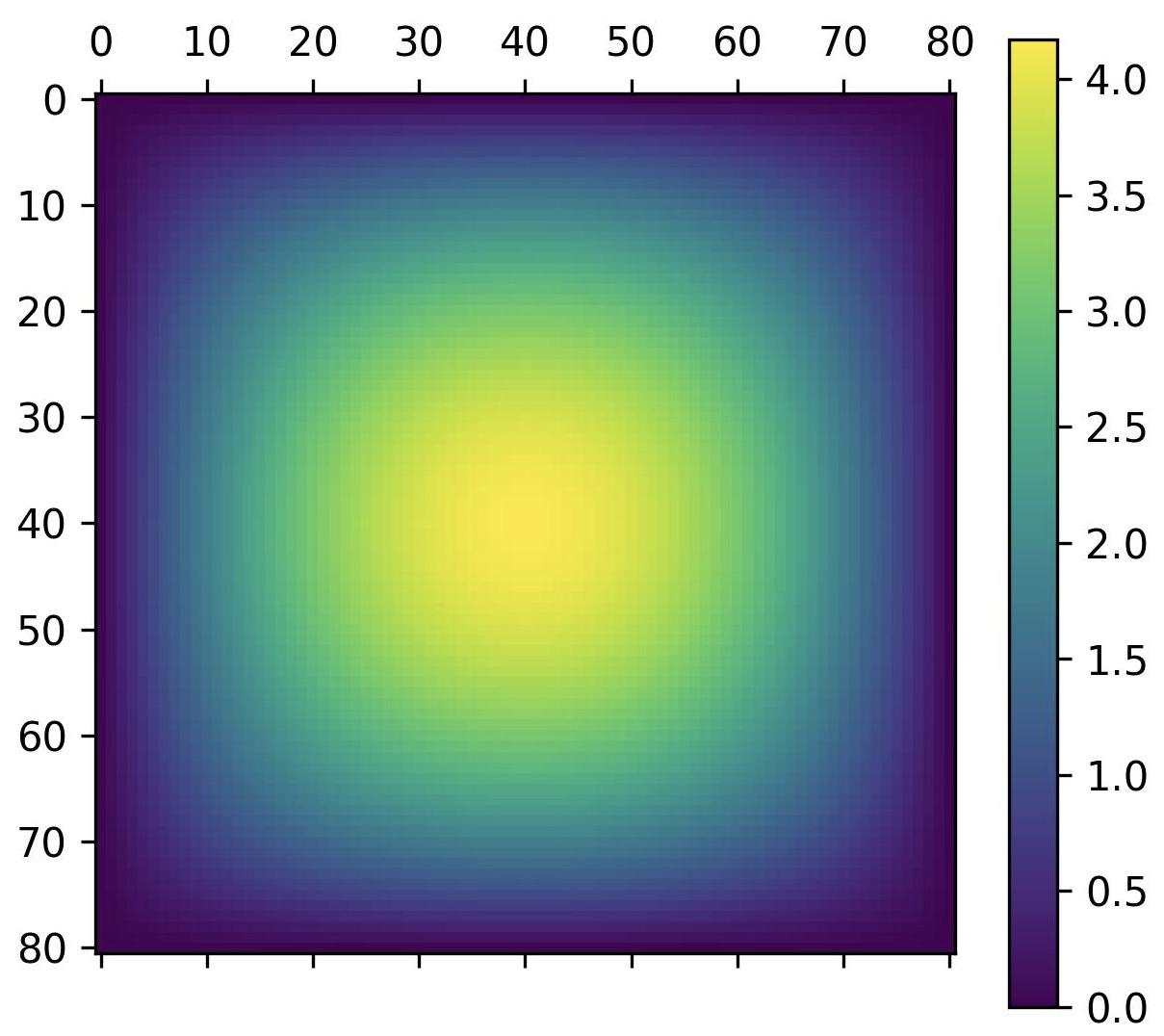} \\
		(c) $\honeoneMC$ 
	\caption{Comparison of h-statistics-based MLMC (top left), classical MLMC (top right) and h-statistics-based MC (bottom) covariance estimation of the stochastic temperature field $\vek{U}^{h_L}(\event)$ on finest level $L=3$}
	\label{Fig:covaccu}
\end{figure}

\begin{figure}[ht!]
	\centering
	\includegraphics[width=0.45\linewidth]{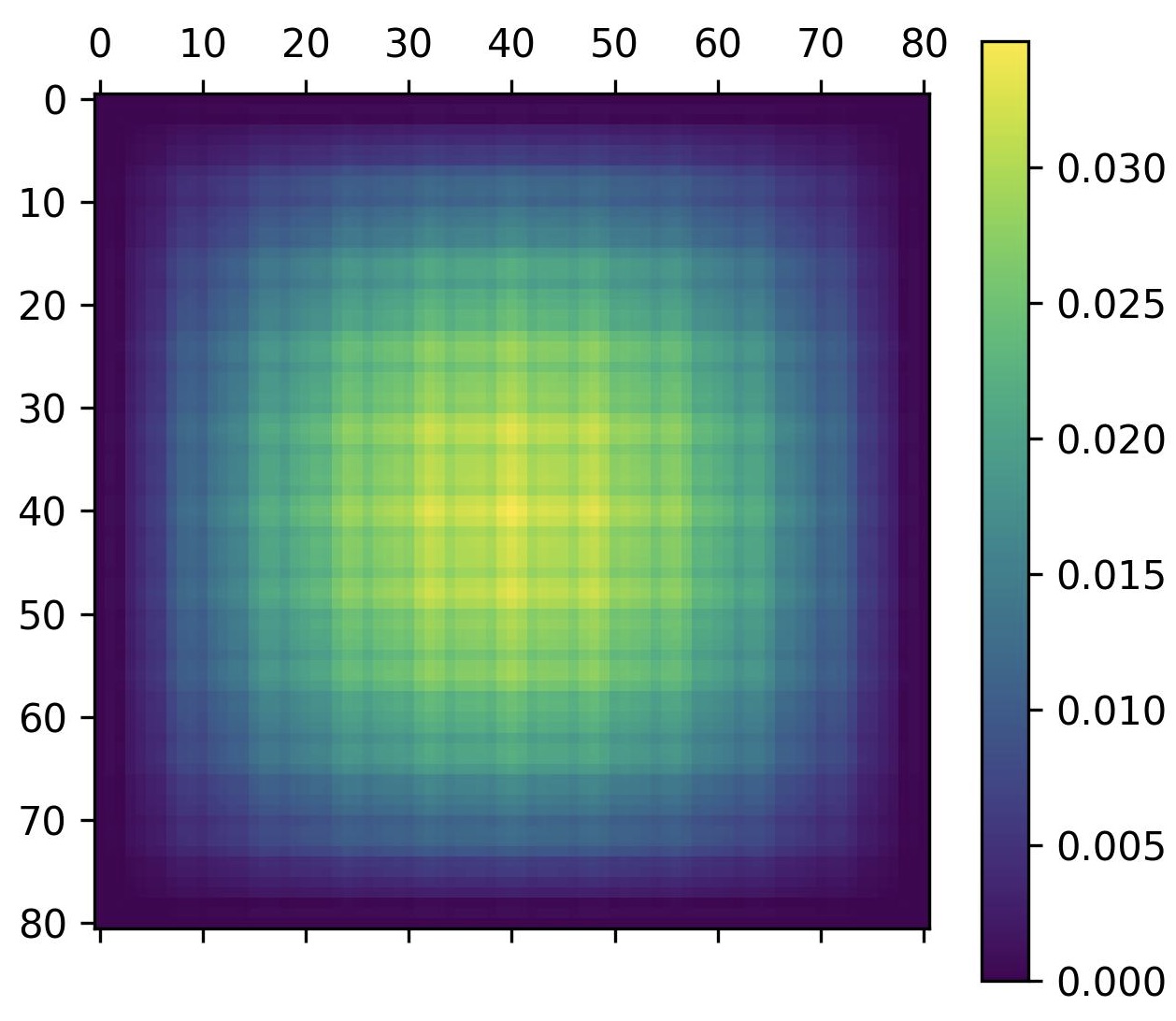} \hfill
	\includegraphics[width=0.45\linewidth]{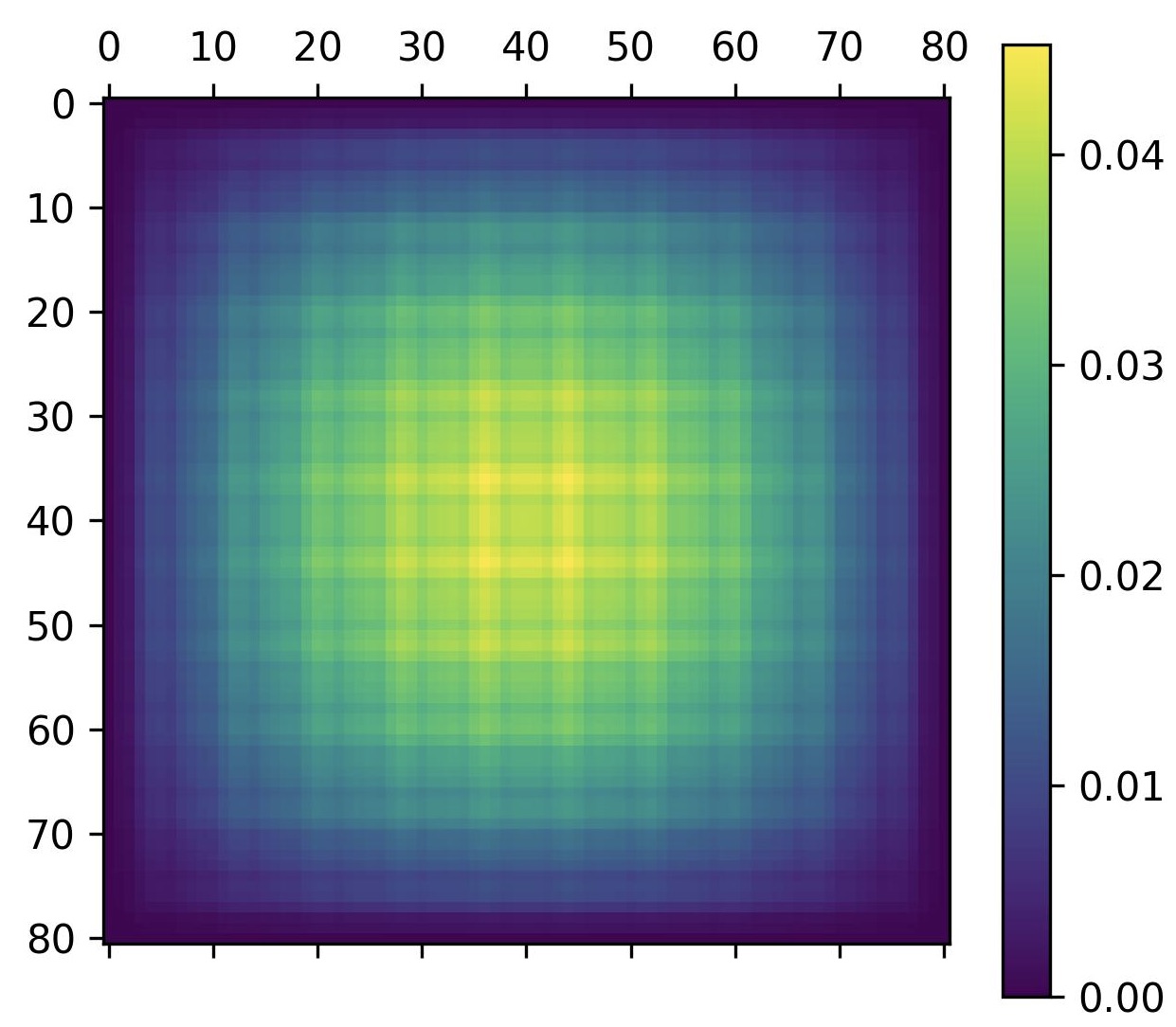} \\
	\rvsi{(a) $\lvert \honeoneML-\honeoneMC \rvert$ \hspace*{190pt}	(b) $\lvert \covML-\honeoneMC \rvert$} \\
	\caption{\rvsi{Comparison of absolute errors (in $\mathrm{K}^4$) between h-statistics-based MLMC and classical MLMC covariance estimators, with respect to h-statistics-based MC estimation}}
	\label{Fig:abserr}
\end{figure}

\ffig{Fig:covaccu} compares the covariance estimates of the temperature for a 1D bar using three different estimators. The results are presented for an absolute stochastic accuracy of $\epssqtwo = 10^{-3}\, \mathrm{K}^4$. In all three sub-figures, a significant covariance in the temperature field is observed at the centre of the spatial domain. It is evident that the accuracy of the statistical covariance estimate is consistent across all estimators.
\rvsi{Furthermore, the absolute errors of both the MLMC covariance estimators (in $\mathrm{K}^4$) is determined with respect to the MC estimate, which is displayed in \ffig{Fig:abserr}. Here, the left plot which signifies the absolute error corresponding to h-statistics based MLMC estimator and the right which is related to the classical MLMC incur small errors. It is also evident that the left plot displays overall smaller error as compared to the right.}

\rvsi{\ftbl{Tbl:relerr} tabulates the total relative error of both the multilevel estimators with regard to MC estimate against all three considered stochastic accuracies. For example, the relative error with respect to h-statistics based MLMC estimator is determined as $\lVert \honeoneML-\honeoneMC \rVert/\lVert\honeoneMC \rVert$. Evidently, both the estimators achieve total relative accuracy of $\le1\%$; more specifically, the estimator $\honeoneML$ achieves better accuracy as compared to the estimator $\covML$ for all stochastic errors except $0.75 \times 10^{-3}$}.
\begin{table}[!ht]
	\centering
	\rvsi{
	\begin{tabular}{llll}
		\specialrule{1pt}{1pt}{1pt}
		 & $10^{-3}$ & $0.75 \times 10^{-3}$  & $0.5 \times 10^{-3}$  \\ \midrule
		$\honeoneML$ & 0.73  & 0.75   &	0.38   \\
		$\covML$ & 1.00  & 0.40  & 0.71  \\
		\specialrule{1pt}{1pt}{1pt}
	\end{tabular}}
	\caption{\rvsi{Comparison of total relative error (in percentage) for each multilevel covariance estimator relative to h-statistics-based MC estimation across three levels of stochastic accuracies (in $\mathrm{K}^4$)}}
	\label{Tbl:relerr}
\end{table}

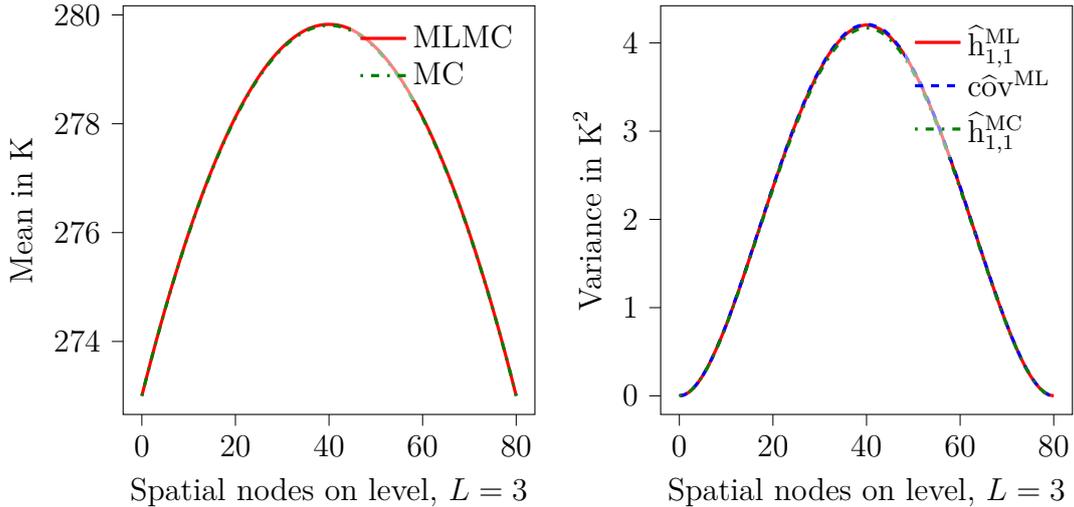
\begin{figure}[!ht]
	\centering
\begin{tikzpicture}
\definecolor{darkgray176}{RGB}{176,176,176}
\definecolor{goldenrod1911910}{RGB}{191,191,0}
\definecolor{green01270}{RGB}{0,127,0}
\definecolor{lightgray204}{RGB}{204,204,204}

\begin{axis}[
width=7cm,
height=7cm,
legend cell align={left},
legend style={anchor=north east, fill opacity=0.5, draw opacity=1, text opacity=1, draw=none},
tick align=outside,
tick pos=left,
x grid style={darkgray176},
xlabel={Spatial nodes on level, $L=3$ },
xmin=-4, xmax=84,
xtick style={color=black},
y grid style={darkgray176},
ylabel={Mean in K},
ymin=272.658653973751, ymax=280.168266551232,
ytick style={color=black},
axis background/.style={fill=white},
]
\addplot [semithick, line width=1.2pt, red]
table {%
0 273
1 273.33796606579
2 273.667015702573
3 273.987653837609
4 274.299375543639
5 274.602076684858
6 274.89586139707
7 275.181234607535
8 275.457691388994
9 275.727388249534
10 275.988168681067
11 276.240537610854
12 276.483990111634
13 276.718422047603
14 276.943937554565
15 277.16104155978
16 277.369229135989
17 277.570656791279
18 277.763168017563
19 277.947267742099
20 278.12245103763
21 278.288613768349
22 278.445860070061
23 278.594694870027
24 278.734613240986
25 278.867771691026
26 278.99201371206
27 279.107844231347
28 279.214758321627
29 279.312651847096
30 279.401628943559
31 279.482194538274
32 279.553843703984
33 279.618732948774
34 279.674705764558
35 279.722267078595
36 279.760911963626
37 279.790536283845
38 279.811244175057
39 279.823540564523
40 279.826920524983
41 279.823540564523
42 279.811244175057
43 279.790536283845
44 279.760911963625
45 279.722267078595
46 279.674705764558
47 279.618732948774
48 279.553843703983
49 279.482194538274
50 279.401628943558
51 279.312651847096
52 279.214758321627
53 279.107844231346
54 278.992013712059
55 278.867771691025
56 278.734613240985
57 278.594694870026
58 278.445860070061
59 278.288613768348
60 278.122451037629
61 277.947267742099
62 277.763168017562
63 277.570656791279
64 277.369229135989
65 277.16104155978
66 276.943937554564
67 276.718422047602
68 276.483990111633
69 276.240537610853
70 275.988168681067
71 275.727388249534
72 275.457691388994
73 275.181234607535
74 274.89586139707
75 274.602076684858
76 274.299375543639
77 273.987653837609
78 273.667015702573
79 273.33796606579
80 273
};
\addlegendentry{MLMC}
\addplot [semithick, line width=1.2pt, green01270, dash pattern=on 1pt off 3pt on 3pt off 3pt]
table {%
0 273
1 273.336178955266
2 273.663847050905
3 273.983004286917
4 274.293650663302
5 274.59578618006
6 274.889410837191
7 275.174524634695
8 275.451127572572
9 275.719219650822
10 275.978800869445
11 276.229871228441
12 276.47243072781
13 276.706479367552
14 276.932017147668
15 277.149044068156
16 277.357560129017
17 277.557565330251
18 277.749059671858
19 277.932043153838
20 278.106515776191
21 278.272477538918
22 278.429928442017
23 278.578868485489
24 278.719297669335
25 278.851215993553
26 278.974623458144
27 279.089520063108
28 279.195905808446
29 279.293780694156
30 279.383144720239
31 279.463997886696
32 279.536340193525
33 279.600171640728
34 279.655492228303
35 279.702301956251
36 279.740600824573
37 279.770388833267
38 279.791665982335
39 279.804432271775
40 279.808687701589
41 279.804432271775
42 279.791665982335
43 279.770388833267
44 279.740600824573
45 279.702301956251
46 279.655492228303
47 279.600171640728
48 279.536340193525
49 279.463997886696
50 279.38314472024
51 279.293780694156
52 279.195905808446
53 279.089520063109
54 278.974623458144
55 278.851215993553
56 278.719297669335
57 278.578868485489
58 278.429928442017
59 278.272477538918
60 278.106515776192
61 277.932043153838
62 277.749059671858
63 277.557565330251
64 277.357560129017
65 277.149044068156
66 276.932017147668
67 276.706479367552
68 276.47243072781
69 276.229871228441
70 275.978800869445
71 275.719219650822
72 275.451127572572
73 275.174524634695
74 274.889410837191
75 274.59578618006
76 274.293650663302
77 273.983004286917
78 273.663847050905
79 273.336178955266
80 273
};
\addlegendentry{MC}
\end{axis}

\end{tikzpicture}
	\centering
\begin{tikzpicture}

\definecolor{darkgray176}{RGB}{176,176,176}
\definecolor{goldenrod1911910}{RGB}{191,191,0}
\definecolor{green01270}{RGB}{0,127,0}
\definecolor{lightgray204}{RGB}{204,204,204}

\begin{axis}[
	width=7cm,
	height=7cm,
legend cell align={left},
legend style={anchor=north east, fill opacity=0.5, draw opacity=1, text opacity=1, draw=none},
tick align=outside,
tick pos=left,
x grid style={darkgray176},
xlabel={Spatial nodes on level, $L=3$ },
xmin=-4, xmax=84,
xtick style={color=black},
y grid style={darkgray176},
ylabel={Variance in $\text{K}^2$},
ymin=-0.210535928727652, ymax=4.42125450328068,
ytick style={color=black}
]
\addplot [semithick, line width=1.2pt, red]
table {%
0 0
1 0.0101662524259399
2 0.0394348075953617
3 0.0871066472441362
4 0.151037394246419
5 0.230084983411698
6 0.322228077534629
7 0.428294356625134
8 0.544911747002904
9 0.66991201796088
10 0.802328251575922
11 0.944362155766529
12 1.09156726033701
13 1.24393207310886
14 1.39900767390107
15 1.56021642801732
16 1.72219019573574
17 1.88270308480734
18 2.04162964321416
19 2.20346020081668
20 2.36205831468738
21 2.51834829081857
22 2.66964839371283
23 2.82136354897689
24 2.96674170359485
25 3.1036277021833
26 3.23260915155583
27 3.35985286428292
28 3.47814456759908
29 3.5891495621814
30 3.69014806000143
31 3.78791547574276
32 3.87492793502362
33 3.94922947266609
34 4.01199227726131
35 4.07044759280862
36 4.1169153565843
37 4.15360713555146
38 4.17795984749868
39 4.19750731155562
40 4.20456590577426
41 4.19750731055729
42 4.17795984749868
43 4.15360713559829
44 4.1169153565843
45 4.07044759291802
46 4.01199227783271
47 3.94922947365747
48 3.87492793502363
49 3.78791547574276
50 3.69014806000143
51 3.58914956366639
52 3.47814456759908
53 3.35985286428292
54 3.23260915155583
55 3.1036277021833
56 2.9667417051565
57 2.82136355022833
58 2.66964839371283
59 2.51834829081857
60 2.36205831468738
61 2.20346020081668
62 2.04162964321416
63 1.88270308480733
64 1.72219019573574
65 1.56021642801732
66 1.39900767390107
67 1.24393207310886
68 1.09156726033701
69 0.944362155922775
70 0.802328249620223
71 0.669912019994074
72 0.544911747002904
73 0.428294356625134
74 0.322228079051384
75 0.230084982375819
76 0.151037394246419
77 0.0871066464849076
78 0.0394348075953539
79 0.0101662524259337
80 0
};
\addlegendentry{$\honeoneML$}
\addplot [semithick, line width=1.2pt, blue, dashed]
table {%
0 0
1 0.0102639109240819
2 0.0398226285433045
3 0.0881386660626042
4 0.153007388190418
5 0.231901370221322
6 0.323661014503167
7 0.429474385700839
8 0.545709134636293
9 0.671585018465343
10 0.805024888190865
11 0.948757492796051
12 1.09788822449787
13 1.24921311853736
14 1.40306846237214
15 1.56355817869596
16 1.724710328076
17 1.88662129634385
18 2.04705166311198
19 2.21129564756256
20 2.37246941936584
21 2.52696405001561
22 2.67634012533094
23 2.82691673253473
24 2.97108303810071
25 3.10966363586862
26 3.24040526257892
27 3.37046689062946
28 3.49167620614109
29 3.60029083199908
30 3.698824136772
31 3.79510958305685
32 3.88059824048814
33 3.95669447836177
34 4.02127958646539
35 4.08267660052528
36 4.13212540447366
37 4.16602173116311
38 4.18756018235288
39 4.20538791935966
40 4.21071857455303
41 4.20538791948155
42 4.18756018235288
43 4.16602173193146
44 4.13212540447366
45 4.08267660054876
46 4.02127958646539
47 3.95669447876891
48 3.88059824048814
49 3.79510958305685
50 3.69882413785791
51 3.60029083199908
52 3.49167620614109
53 3.37046689183473
54 3.24040526257892
55 3.10966363598108
56 2.9710830381007
57 2.82691673174735
58 2.67634012533094
59 2.52696405001561
60 2.37246941936584
61 2.21129564756256
62 2.04705166311198
63 1.88662129634385
64 1.724710328076
65 1.56355817869596
66 1.40306846247684
67 1.24921311373954
68 1.09788822391552
69 0.948757492269826
70 0.805024887819894
71 0.671585019134454
72 0.545709134636292
73 0.429474385700839
74 0.323661014910488
75 0.231901370221322
76 0.153007388190418
77 0.088138665760505
78 0.0398226285432996
79 0.0102639109240717
80 0
};
\addlegendentry{$\covML$}
\addplot [semithick, line width=1.2pt, green01270, dash pattern=on 1pt off 3pt on 3pt off 3pt]
table {%
0 0
1 0.0101657688500228
2 0.0396401459276003
3 0.0869180566585588
4 0.150533519314806
5 0.229059645014332
6 0.32110863772121
7 0.425331794245597
8 0.54041950424373
9 0.665101250217929
10 0.798145607516599
11 0.938360244334223
12 1.08459192171137
13 1.23572649353469
14 1.39068890653692
15 1.54844320029687
16 1.70799250723944
17 1.8683790526356
18 2.02868415460243
19 2.18802822410307
20 2.34557076494674
21 2.50051037378875
22 2.6520847401305
23 2.79957064631946
24 2.94228396754919
25 3.07957967185932
26 3.21085182013559
27 3.33553356610979
28 3.45309715635981
29 3.56305393030962
30 3.66495432022928
31 3.75838785123491
32 3.84298314128873
33 3.91840790119905
34 3.98436893462024
35 4.04061213805277
36 4.08692250084318
37 4.12312410518411
38 4.14908012611427
39 4.16469283151844
40 4.16990358212751
41 4.16469283151844
42 4.14908012611425
43 4.1231241051841
44 4.08692250084316
45 4.04061213805274
46 3.98436893462021
47 3.91840790119901
48 3.84298314128869
49 3.75838785123487
50 3.66495432022923
51 3.56305393030958
52 3.45309715635976
53 3.33553356610974
54 3.21085182013554
55 3.07957967185928
56 2.94228396754914
57 2.79957064631942
58 2.65208474013046
59 2.50051037378871
60 2.3455707649467
61 2.18802822410303
62 2.0286841546024
63 1.86837905263557
64 1.70799250723941
65 1.54844320029684
66 1.39068890653689
67 1.23572649353467
68 1.08459192171135
69 0.938360244334204
70 0.798145607516581
71 0.665101250217914
72 0.540419504243717
73 0.425331794245585
74 0.321108637721201
75 0.229059645014325
76 0.1505335193148
77 0.0869180566585555
78 0.0396401459275986
79 0.0101657688500223
80 0
};
\addlegendentry{$\honeoneMC$}
\end{axis}

\end{tikzpicture}
	\caption{Comparison of classical MLMC and MC mean (on the left) and the comparison of all three covariance estimators for variance estimation (on the right) of the stochastic temperature field $\vek{U}^{h_L}(\event)$ on finest level $L=3$}
	\label{Fig:meanvaraccu}
\end{figure}

For further analysis, we estimate the MLMC and MC mean of the temperature profile, $u^h({x},\event)$, with a sampling error of $\epssqtwo = 10^{-3}\, \mathrm{K}^2$ using a conventional approach \cite{Cliffe2011, giles2008}.
We calculate the diagonal of the resulting covariance matrices from \ffig{Fig:covaccu} to determine the variance of the stochastic response. The corresponding results are presented in \ffig{Fig:meanvaraccu}.
In the left plot, we observe a maximum mean value of approximately 279 K at the center of the 1D geometry, where both MLMC and MC estimates are nearly identical. \rvsi{The total relative error of MLMC mean with regard to MC mean is 0.002$\%$.}
Conversely, the right plot illustrates the influence of input material uncertainty, $\kappa(\event)$, on the temperature, $u^h(x,\event)$. The maximum variance is also observed in the same region as the maximum mean value. Additionally, all three estimators yield similar results; \rvsi{specifically, the relative difference of $\honeoneML$  and $\covML$ with respect to $\honeoneMC$ are noted in \ftbl{Tbl:relerr_var}. Similar to the performance of covariance estimators, the corresponding variance estimators achieve a total relative accuracy of $\le1\%$ for different sampling accuracies.}
\begin{table}[!ht]
	\centering
	\rvsi{
	\begin{tabular}{llll}
		\specialrule{1pt}{1pt}{1pt}
		& $10^{-3}$ & $0.75 \times 10^{-3}$  & $0.5 \times 10^{-3}$  \\ \midrule
		$\honeoneML$ & 0.74  & 0.74   &	0.38  \\
		$\covML$ & 1.0  & 0.40  & 0.71  \\
		\specialrule{1pt}{1pt}{1pt}
	\end{tabular}}
	\caption{\rvsi{Comparison of total relative error (in percentage) in variance estimation by each multilevel covariance estimator relative to h-statistics-based MC estimation across three levels of stochastic accuracies (in $\mathrm{K}^4$)}}
	\label{Tbl:relerr_var}
\end{table}

Note that the specified sampling accuracies for mean and covariance/variance are of similar magnitudes but belong to two different scales, namely $\mathrm{K}^2$ and $\mathrm{K}^4$. Therefore, a direct comparison of accuracies between these statistics is not straightforward in this scenario, given the difference in scales. Since the focus of this study centers on the performance analysis of covariance estimators rather than a cross-comparison of different statistics, we direct readers to \cite{shivanand_scale-invariant_2023} for a description of normalized error estimates in MLMC and MC mean and variance estimators, which addresses this issue.

\rvsi{Finally, the summary of the given and estimated maximum values of sampling errors of MLMC covariance estimators (in $\mathrm{K}^4$) is tabulated in \ftbl{Tbl:achaccu}. It can be observed that the achieved maximum stochastic accuracies fall within the prescribed error limits of $\epssqtwo$.}
\begin{table}[!ht]
	\centering
	\rvsi{
	\begin{tabular}{llll}
		\specialrule{1pt}{1pt}{1pt}
		& {$10^{-3}$} & $0.75 \times 10^{-3}$  & $0.5 \times 10^{-3}$  \\ \midrule
		$\honeoneML$ & $0.95 \times 10^{-3}$  &  $0.72 \times 10^{-3}$  &	$0.48 \times 10^{-3}$   \\
		$\covML$ & $0.98 \times 10^{-3}$  &  $0.73 \times 10^{-3}$ & $0.49 \times 10^{-3}$  \\
		$\covMC$ & $0.75 \times 10^{-3}$  &  $0.56 \times 10^{-3}$ & $0.36 \times 10^{-3}$  \\
		\specialrule{1pt}{1pt}{1pt}
	\end{tabular}}
	\caption{\rvsi{Comparison of estimated sampling accuracies among different covariance estimators across three levels of stochastic errors (in $\mathrm{K}^4$)}}
	\label{Tbl:achaccu}
\end{table}

\subsubsection{Efficiency}

\ftbl{Tbl:costsingle} provides an overview of the computational time required to obtain a single deterministic output sample on each mesh level $l$. 
\begin{table}[!ht]
	\centering
	\begin{tabular}{llll}
		\specialrule{1pt}{1pt}{1pt}
		Mesh level & {\begin{tabular}[c]{@{}l@{}}$\honeoneML$,  \\ $\mathcal{C}(\Zl)$ \end{tabular}}& {\begin{tabular}[c]{@{}l@{}}$\covML$,  \\ $\mathcal{C}(\Yl)$ \end{tabular}} & {\begin{tabular}[c]{@{}l@{}}$\honeoneMC$,  \\ $\mathcal{C}(u^{h_L})$ \end{tabular}}  \\ \midrule
		0& 0.15  & 0.15   &	-    \\
		1& 1.11   & 1.11 & -    \\
		2& 1.96  &1.96  & -    \\
		3&2.94 & 2.94 & 1.94  \\
		\specialrule{1pt}{1pt}{1pt}
	\end{tabular}
	\caption{Computational time to run one sample in millisecond for different estimators on each mesh level $l$}
	\label{Tbl:costsingle}
\end{table} 
To determine the maximum number of samples, denoted as $\max(N_l)$, for the h-statistics-based MLMC covariance estimator, $\honeoneML$, we utilize the expression given in \feq{Eq:tau} within \feq{eq:Nl}:
\begin{equation}
	\max(N_l) = \max \left( \frac{2}{\epssq} \left( \sum_{\meshlevel=0}^{\numofmesh} \left({{\vech(\vloneone)\mathcal{C}(\Zl)}}\right) ^{\circ \frac{1}{2}}\right) \circ \left({\frac{\vech(\vloneone)}{\mathcal{C}(\Zl)}}\right)^{\circ\frac{1}{2}}\right).
\end{equation}
Here, the symbol $\circ$ represents the Hadamard or element-wise operator. A similar procedure is followed to determine $\max(N_l)$ for the conventional estimator, $\covML$. Additionally, the maximum number of samples generated at $L=3$ for the MC estimator, $\honeoneMC$, is determined as:
\begin{equation}
	\max(\samplesize_L) = \max\left(\frac{\vech(\voneoneMC)}{\epssqtwo}\right).
\end{equation}
The corresponding results, considering varying sampling accuracies, are presented in \ftbl{Tbl:nmax}. 
\begin{table}[!ht]
	\centering
		\begin{tabular}{llll}
			\specialrule{1pt}{1pt}{1pt}
			Mesh level & $\honeoneML$ & $\covML$  & $\honeoneMC$  \\ \midrule
			0& 66822   & 82417    &	-    \\
			1& 397   & 595  & -    \\
			2& 82  &95   & -    \\
			3&50 & 50 & 84582   \\
			\specialrule{1pt}{1pt}{1pt}
		\end{tabular} \hspace{0.5cm} 
		\begin{tabular}{llll}
			\specialrule{1pt}{1pt}{1pt}
			Mesh level & $\honeoneML$ & $\covML$  & $\honeoneMC$  \\ \midrule
			0& 86879  &  110072  &	-    \\
			1&  600  &  761& -    \\
			2& 101 &  117& -    \\
			3& 50 &  50& 112777 \\
			\specialrule{1pt}{1pt}{1pt}
		\end{tabular} \\
		(a) $\epssqtwo = 10^{-3}\, \mathrm{K}^4$ \hspace*{80pt} (b) $\epssqtwo = 0.75 \times 10^{-3}\, \mathrm{K}^4$ \\ 
		\begin{tabular}{llll}
			\specialrule{1pt}{1pt}{1pt}
			Mesh level & $\honeoneML$ & $\covML$  & $\honeoneMC$  \\ \midrule
			0& 129008  &  164750  &	-    \\
			1&  849 & 1132 & -    \\
			2& 179  & 229& -    \\
			3& 50 &53 &  169165 \\
			\specialrule{1pt}{1pt}{1pt}
		\end{tabular} \\
		(c) $\epssqtwo = 0.5 \times 10^{-3}\, \mathrm{K}^4$ 
	\caption{Maximum number of MC samples on each mesh level $l$ for different estimators with respect to varying stochastic accuracies}
	\label{Tbl:nmax}
\end{table}
It is worth noting that $\max(N_l)$ in MLMC estimators exhibits a monotonic decrease as $l$ increases. Most importantly, for all the explored sampling accuracies, it is evident that the estimator $\honeoneML$ requires fewer samples compared to $\covML$.

\begin{figure}[!ht]
	\centering
\begin{tikzpicture}

\definecolor{darkgray176}{RGB}{176,176,176}
\definecolor{steelblue31119180}{RGB}{31,119,180}
\definecolor{darkorange25512714}{RGB}{255,127,14}
\definecolor{lightgray204}{RGB}{204,204,204}
\definecolor{green01270}{RGB}{0,127,0}

\begin{axis}[
width=7cm,
height=7cm,
legend cell align={left},
legend style={at={(1.05,1.05)}, anchor=north west, legend cell align=left, align=left, draw=none},
tick align=outside,
tick pos=left,
x grid style={darkgray176},
xlabel={Type of estimator},
xtick style={color=black},
xtick={0,1,2},
xticklabels={$\honeoneML$,$\covML$,$\honeoneMC$},
y grid style={darkgray176},
ylabel={Cost in second},
ymin=1.1774652532615, ymax=330,
ytick style={color=black}
]
\addplot [semithick, line width=1.2pt, steelblue31119180,dash pattern=on 1pt off 3pt on 3pt off 3pt, mark=square, mark size=3, mark options={solid}]
table {%
0 10.688540580424
1 13.2534232071671
2 164.21891772144
};
\addlegendentry{$\epssqtwo = 10^{-3}\, \mathrm{K}^4$}

\addplot [semithick, line width=1.2pt, green01270, dashed, mark=square, mark size=3, mark options={solid}]
table {%
	0 13.9345321999677
	1 17.5946890831855
	2 218.960498496972
};
\addlegendentry{$\epssqtwo = 0.75\times10^{-3}\, \mathrm{K}^4$}

\addplot [semithick, line width=1.2pt, darkorange25512714, mark=square, mark size=3, mark options={solid}]
table {%
	0 20.630665236779
	1 26.3682672435464
	2 328.439776977932
};
\addlegendentry{$\epssqtwo = 0.5\times10^{-3}\, \mathrm{K}^4$}
\end{axis}

\end{tikzpicture}
	\centering
	\caption{Comparison of total computational cost to estimate sample covariance on mesh $L=3$ for different estimators}
	\label{Fig:cost}
\end{figure}
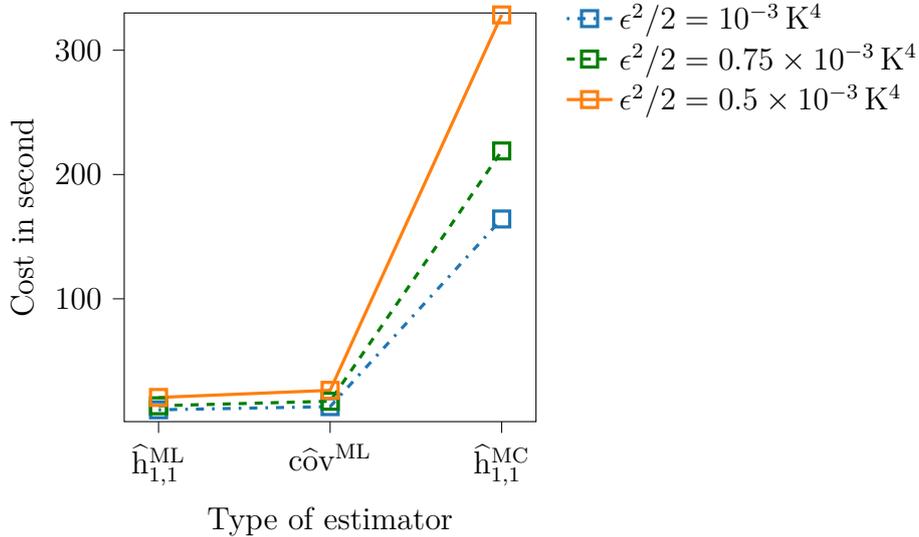

To determine the total computational cost of all three estimators, we formulate it as follows:
\begin{align}
	\mathcal{C}(\honeoneML) &= \sum_{\meshlevel=0}^{\LL} \max(\Nl) \mathcal{C}(\Zl), \nonumber \\
	\mathcal{C}(\covML) &= \sum_{\meshlevel=0}^{\LL} \max(\Nl) \mathcal{C}(\Yl), \\
	\mathcal{C}(\honeoneMC) &= \max(\samplesize_L) \mathcal{C}(u^{h_L}). \nonumber
\end{align}
In this context, \ffig{Fig:cost} provides an overview of the total computational cost for all three covariance estimators concerning the given stochastic accuracies.
Clearly, both MLMC estimates exhibit a faster convergence rate compared to the MC approach. Furthermore, the computational cost of the h-statistics-based $\honeoneML$ estimator is notably lower than that of the conventional $\covML$ estimator.
Specifically, the h-statistics-based $\honeoneML$ estimator is $\approx16\times$ faster than the h-statistics-based $\honeoneMC$ estimator for all the studied accuracies, resulting in cost savings of $\approx94\%$.
Similarly, the closed-form and unbiased sampling error-based $\honeoneML$ is  $\approx 1.3\times$ faster than the conventional upper bound and biased sampling error-based $\covML$, leading to cost savings of  $\approx22\%$ when choosing $\honeoneML$ over $\covML$.

\section{Conclusion}\label{Conclusion}

A novel approach to estimate the unbiased covariance of random variables, \rvsn{with minimal variance}, using both Monte Carlo (MC) and Multilevel Monte Carlo (MLMC) methods with h-statistics is introduced. Traditional MC covariance estimation relies on a closed-form, yet biased, sampling error; in contrast, MLMC covariance estimators depend on biased and upper-bounded stochastic accuracy. To enhance the overall accuracy and efficiency of covariance estimation, we derive unbiased and closed-form versions of the mean square error estimates using h-statistics within both the MC and MLMC algorithms. \rvsn{This leads to potentially achieving sharper error estimates as compared to the error bounds of traditional estimators, particularly in MLMC estimation.}

The proposed h-statistics-based procedures, along with the conventional MLMC method, are tested by applying them to Poisson's equation in the form of a simple 1D steady-state heat equation. The analysis incorporates uncertainty in the material properties, specifically the thermal conductivity, modelled as a log-normal random variable. These methods propagate the uncertainty to estimate the covariance, as well as statistics, such as the mean and variance of the temperature field.

Finally, we compare the computational efficiencies of the MLMC and MC covariance estimators for different stochastic accuracies. The h-statistics-based MLMC estimator significantly outperformed the corresponding MC method in terms of the computational cost. It is also evident that the h-statistics-based MLMC method required fewer samples than the conventional MLMC, highlighting the computational advantages of the unbiased and closed-form sampling error-based MLMC covariance estimator over biased and upper-bounded MLMC estimates. Furthermore, the difference in accuracy between the MLMC and MC methods is small.

\section*{Acknowledgments}
This work was supported by Wave 1 of The UKRI Strategic Priorities Fund under the EPSRC Grant EP/W006022/1, particularly the Digital Twins: Complex Systems Engineering theme within that grant \& The Alan Turing Institute. The author extends their sincere gratitude to Prof. Fehmi Cirak and the Department of Engineering at the University of Cambridge for hosting and facilitating initial discussions on the deterministic solver and providing software support.

\bibliographystyle{abbrv}
\bibliography{MLMC_cov}


\end{document}